\documentclass[xypic,amscd,syntonly,amssymb,verbatim,12pt]{amsart}

\usepackage{graphicx}
\usepackage[all]{xy}
\usepackage{amssymb}
\usepackage{amsmath}
\usepackage[mathscr]{euscript}

\usepackage{epsfig}
\theoremstyle{plain}
\newtheorem{Thm}{Theorem}
\newtheorem{Prop}{Proposition}
\newtheorem{Cor}{Corollary}
\newtheorem{Lem}{Lemma}
\newtheorem{Def}{Definition}

 \theoremstyle{definition}
\theoremstyle{remark}

\numberwithin{equation}{section}

\begin{document}
  \title {Gauge fields on coherent sheaves}

 \author{ Andr\'{e}s   Vi\~{n}a}
\address{Departamento de F\'{i}sica. Universidad de Oviedo.    Garc\'{\i}a Lorca 18.
     33007 Oviedo. Spain. }
\email{vina@uniovi.es}
  \keywords{Yang-Mills fields, holomorphic connections, coherent sheaves, flat gauge fields}

 \maketitle
\begin{abstract}
  Given a flat gauge field $\nabla$ on a vector bundle $F$ over a manifold $M$ we deduce a necessary and sufficient condition for the field $\nabla+ E$, with $E$ an ${\rm End}(F)$-valued $1$-form, to be a Yang-Mills field. For each curve of Yang-Mills fields on $F$ starting at $\nabla$, we define a cohomology class of $H^2(M,\,{\mathscr P})$, with ${\mathscr P}$ the sheaf of $\nabla$-parallel sections of $F$. This cohomology class vanishes when the curve consists of flat fields.  We prove the existence of a curve of Yang-Mills fields on a bundle over the torus $T^2$ connecting two vacuum states.

We define holomorphic and meromorphic gauge fields on a coherent sheaf and the corresponding Yang-Mills functional. In this setting, we analyze the Aharonov-Bohm effect and the Wong equation.

\end{abstract}
   \smallskip
 MSC 2010: 53C05, 14F06, 58E15, 32C35

\section {Introduction} \label{S:intro}
 
 From the mathematical point of view, a gauge field 
is a connection on a coherent sheaf, and the field strength  is 
the curvature of the connection   \cite{C-J, Hamilton, Naber}. 

The set of connections on a coherent sheaf ${\mathscr F}$ over the manifold $M$ is an affine space associated to the vector space 
${\rm Hom}({\mathscr F},\,{\mathscr A}^1\otimes{\mathscr F})$, where ${\mathscr A}^1$ is the sheaf of $1$-forms on $M$. On the other hand, each element 
$E\in \Gamma(M,\,{\mathscr A}^1\otimes{\mathscr End}({\mathscr F}))$ determines a vector   ${\mathcal E}\in{\rm Hom}({\mathscr F},\, {\mathscr A}^1\otimes{\mathscr F})$. But this correspondence is not necessarily $1$ to $1$. 
The map $E\mapsto {\mathcal E}$ is bijective, if  ${\mathscr F}$ is locally free  (Lemma \ref{L:1}) or when ${\mathscr F}$ is a reflexive sheaf (Lemma \ref{Lem:finite}) (the reflexive sheaves might be thought as ``vector bundles with singularities'' \cite[page 121]{Hartshorne1}).
We prove some properties of the connections on coherent sheaves, as the relation between the curvatures of two connections
 (Proposition \ref{KHatNabla}). By means of the identification $E\leftrightarrow{\mathcal E}$, these properties give rise in the setting of locally free sheaves 
to well-known formulas for the connections on vector bundles. However, we will rewrite these formulas in a  suitable way
 to prove some results   mentioned below. 

When the coherent sheaf  is   a Hermitian vector bundle $F$ over a Riemannian manifold $M$, the metrics on $F$ and $M$ determine a norm for the ${\rm End}(F)$-valued $2$-forms on $M$. In this situation, one defines the Yang-Mills functional,  $\sf{YM}$, which associates 
to  each gauge on $F$ compatible with the Hermitian metric the norm of its strength.

 Some vector bundles admit   flat gauge fields; i. e. connections whose curvature vanishes.
In this case the Yang-Mills functional   takes its minimum value at those fields, obviously. They are the vacuum states of the theory \cite[p. 447]{Hamilton}.
  The Yang-Mills fields are the stationary ``points'' of $\sf{YM}$. The calculus of variations shows that those fields must satisfy the Yang-Mills equation (\ref{YMequation}), which is the corresponding motion equation.
	  
   Given a gauge field $\nabla$ on $F$, any other connection on $F$ can be written as $\nabla+E$, where $E$ is  an ${\rm End}(F)$-valued  $1$-form on $M$, as we said. Assumed that $\nabla$ is a flat gauge field,   we give, in Theorem \ref{Th:YMHat}, a necessary and sufficient condition for $\nabla+E$ to be a Yang-Mills field.  

The curves in the space of Yang-Mills fields, with initial point at a flat configuration, are possible evolutions of the field from a vacuum state. We  consider a curve $\nabla_t$ of Yang-Mills fields, with $\nabla_0$ a flat connection. We prove that such a curve determines an element of the   cohomology $H^2(M,\,{\mathscr P})$, where ${\mathscr P}$ is the sheaf of  sections of $F$ parallel with respect to the flat connection $\nabla_0$.  If the curve $\nabla_t$ consists of flat gauge fields, then the corresponding cohomology class vanishes (Proposition \ref{CEflat}).

On the trivial bundle $T\times {\mathbb C}^2$ over the $2$-torus $T$, we construct a curve $\{\nabla_t\}_{t\in[0,\,1]}$ of Yang-Mills fields, such that the only flat connections are the corresponding to $t=0$ and $t=1$ (Proposition \ref{curveconn}). This curve connects two minimums of the Yang-Mills functional by means of stationary points of this functional.
In physical terms, the system could   evolve    from  a vacuum state  to another one,  through a ``continuous sequence'' of non-vacuum configurations.


As it is well-known a gauge field on a smooth vector bundle allows us to define the parallel transport along smooth curves. In the complex analytic setting, a gauge field should determine a parallel transport along curves  whose tangent field is holomorphic. This involves that the connection must satisfy the  Leibniz's rule for the product of holomorphic functions and sections. That is, the connection must be holomorphic. On the other hand,
 the branes of type $B$ on a complex analytic manifold $X$ are the objects of $ D^b(X)$, the bounded derived category of coherent sheaves
 \cite[Sect. 5.4]{Aspin} \cite[Sect. 5.3]{Aspin-et}.
 Thus,  the coherent sheaves are particular membranes and therefore the study of gauge fields on coherent sheaves  is of interest from the point of view of theoretical Physics.
 

A holomorphic gauge field on a coherent sheaf ${\mathscr F}$ is a splitting of the extension of ${\mathscr F}$ defined by the $1$-jet sheaf 
 ${\mathscr J}({\mathscr F})$ of ${\mathscr F}$.  When ${\mathscr F}$ is a locally free sheaf, the existence of such a splitting  is equivalent to the vanishing of the Atiyah class of the sheaf \cite{Atiyah}.
 For example, the skyscraper sheaves do not admit holomorphic gauge fields  (Proposition \ref{P:Skyscraper}). By contrast,  as a consequence of Cartan's B theorem, we deduce that the holomorphic bundles on a Stein manifold admit holomorphic gauge fields (Proposition \ref{Stein}).

 As the cohomology groups of coherent sheaves over a compact complex  manifold $X$ have finite dimension, the holomorphic gauge fields on a holomorphic vector bundle (if they exist) form a finite dimensional affine space. Using this fact, we prove that  if the Hodge numbers of $X$ satisfies $h^{2,0}(X)\leq h^{1,0}(X)$, then any ``generic''  holomorphic Yang-Mills field on a Hermitian line bundle is flat  (Corollary \ref{Hodge}).
 
Some geometric phases \cite{Bohm_et_al} can be considered in the frame of the meromorphic gauge fields on coherent sheaves.	
In the context of meromorphic connections, we study the Aharonov-Bohm effect \cite{Bohm_et_al, C-J} and the Wong equation \cite{Bala}. 
The Aharonov-Bohm effect can be analyzed 
 by considering  a meromorphic gauge field on a trivial vector bundle over the manifold $X$, on which the electron is moving. More precisely, a meromorphic flat connection with a logarithmic pole at a point $*\in X$.  By means of the connection, one constructs the respective de Rham complex \cite{Hotta}, and the Riemann-Hilbert correspondence (equivalently, the parallel transport) determines a representation of $\pi_1\big({X}\setminus\{*\}\big)$. The phase factor  acquired by the wave function of the electron,  once it has circumvented the solenoid, is determined by that  representation.  
In this setting, the necessary conditions for the respective shift phase to occur clearly appear.

  A   similar study will be done for the variation of the spin of a particle that is moving in presence of a gauge field. For this purpose, we will transform the Wong equation, which governs the variation of spin, in a parallel transport equation.
  According to the gauge field singularities and the  topological properties of the manifold on which the particle is moving, we will deduce the value of different phase factors.


The article is organized as follows. In the first part of Section \ref{S:YMandF}, we introduce notations and deduce  some properties    of the curvature of connections on coherent sheaves to be used below.
   In Subsection \ref{Ss:Yang-Mills_fields}, the above mentioned Theorem \ref{Th:YMHat} is proved. Subsection \ref{Ss:curves_Y-M} is concerned with curves of Yang-Mills fields starting from a vacuum state, and we define the elements of $H^2(X,\,{\mathscr P})$ which they determine. In Subsection \ref{Ss:Connecting}, one constructs a curve of non-flat Yang-Mills fields on the bundle $T\times {\mathbb C}^2\to T$ that  connects two different vacuum states.

In Subsection \ref{Jet1}, the construction of the first jet bundle $J^1(F)$ of a holomorphic bundle $F$ is revised. We show that the existence of a holomorphic connection on $F$ is equivalent to the vanishing in the group ${\rm Ext}^1(F,\, \Omega^1(F))$ of the extension 
$0\to\Omega^1(F)\to J^1(F)\to F \to 0,$ with $\Omega^1(F)$ the sheaf of $F$-valued holomorphic 1-forms.
We prove also Proposition \ref{Stein} and Corollary \ref{Hodge}, above mentioned.  
In Subsection  \ref{connections_on:Coherent}, we generalize, following \cite{Deligne}, the definition of holomorphic connection to coherent sheaves, and then define the functional Yang-Mills on  these fields in Subsection \ref{Ss:HoloYM}.
The meromorphic fields are introduced in Subsection \ref{meromorphicgauge}, 
 and in that context   the Aharonov-Bohm effect and the Wong equation are analyzed.


\section{Yang-Mills fields and flat fields}\label{S:YMandF}

By  $M$ we denote a compact, connected manifold.  We set ${\mathscr A}^k{}$ for the sheaf of differential $k$-forms on $M$,
and ${\mathscr A}{}$ for the sheaf of differentiable functions. 
 
 Let ${\mathscr F}$ be a coherent  ${\mathscr A}{}$-module. 
A  connection on ${\mathscr F}$ is a  morphism of abelian sheaves
\begin{equation}\label{Def:Connection}
\nabla:{\mathscr F}\to  {\mathscr A}^1 \otimes_{{\mathscr A}{}} {\mathscr F},
\end{equation}
satisfying the Leibniz's rule
\begin{equation}\label{Leibniz}
\nabla(g\sigma)=dg\otimes \sigma+g\nabla \sigma,
\end{equation}
where $g$ and $\sigma$ are sections of ${\mathscr A}{}$ and ${\mathscr F}$ (resp.) defined on an open subset of $M.$

The connection $\nabla$ can be extended to a morphism
$$\nabla^k: {\mathscr A}^{k}  \otimes_{{\mathscr A}{}} {\mathscr F} \to {\mathscr A}^{k+1} \otimes_{{\mathscr A}{}} {\mathscr F} ,$$
in the usual way
\begin{equation}\label{nabla_k}
\nabla^k(\beta\otimes\sigma)=   d\beta \otimes \sigma +(-1)^k\beta\wedge\nabla\sigma,   
\end{equation}
 where $\beta\wedge\nabla\sigma$ is the image of  $\beta\otimes\nabla\sigma$ by the natural morphism 
 \begin{equation}\label{natural}
 {\mathcal N}:{\mathscr A}^k{}\otimes_{{\mathscr A}{}}\otimes\big({\mathscr A}^1{}\otimes_{{\mathscr A}{}}{\mathscr F}\big)\to {\mathscr A}^{k+1}{}\otimes_{{\mathscr A}{}}{\mathscr F}.
 \end{equation}

The curvature of $\nabla$ is the morphism of ${\mathscr A}{}$-modules
$${\mathcal K}_{\nabla}:=\nabla^1\circ\nabla:{\mathscr F}\to {\mathscr A}^2 \otimes_{{\mathscr A}{}}{\mathscr F}.$$
It is well known 
that
\begin{equation}\label{integrable}
\nabla^{k+1}\circ\nabla^{k}(\alpha\otimes\sigma)=\alpha \wedge {\mathcal K}_{\nabla}(\sigma).
 \end{equation}
The connection $\nabla$ is called flat (or integrable) if ${\mathcal K}_{\nabla}=0$.

\smallskip

 If $\nabla$ and $\tilde\nabla$ are connections on ${\mathscr F}$, then 
 $$\Hat\nabla-\nabla\in{\rm Hom}_{{\mathscr A}{}}\big({\mathscr F},\,{\mathscr A}{}^1\otimes_{{\mathscr A}{}}{\mathscr F
 }  \big)= \Gamma\big(M,\, {\mathscr Hom}_{{\mathscr A}{}}\big({\mathscr F},\,{\mathscr A}{}^1\otimes_{{\mathscr A}{}}{\mathscr F}  \big) \big),$$
 where ${\mathscr Hom}_{{\mathscr A}} (\,.\,,\,.\,)$ denotes the respective sheaf of germs of ${\mathscr A}$-morphisms \cite[p 87]{Kas-Sch}.
 Conversely, given  
 ${\mathcal E}\in{\rm Hom}_{{\mathscr A}{}}\big({\mathscr F},\,{\mathscr A}{}^1\otimes_{{\mathscr A}{}}{\mathscr F}  \big)$
 and a connection $\nabla$ on ${\mathscr F}$ , then $\Hat\nabla:=\nabla+{\mathcal E}$ is a connection on ${\mathscr F}$. 
On the other hand, ${\mathcal E}$ induces a natural map
\begin{equation}\label{mathcalBk}
{\mathcal E}^k:{\mathscr A}^k{}\otimes_{{\mathscr A}{}}{\mathscr F}\to {\mathscr A}^{k+1}{}\otimes_{{\mathscr A}{}}{\mathscr F},
\;\;\; \alpha\otimes\sigma\mapsto {\mathcal N}\big(\alpha\otimes (-1)^k{\mathcal E}(\sigma)\big).
\end{equation}

Hence, given $\nabla$ and ${\mathcal E}$ one has two morphisms of abelian sheaves $\nabla^1,\,{\mathcal E}^1:{\mathscr A}^1{}\otimes_{{\mathscr A}{}}{\mathscr F}\to {\mathscr A}^{2}{}\otimes_{{\mathscr A}{}}{\mathscr F},$
and  
  five  morphisms ${\mathscr F}\to{\mathscr A}{}^2\otimes_{{\mathscr A}{}}{\mathscr F}$,  
 \begin{equation}\label{nablaB}
{\mathcal K}_{\nabla},\; \;\;{\mathcal K}_{\Hat\nabla},\;\;\; \nabla^1\circ{\mathcal E},\;\;\;
   {\mathcal E}^1\circ\nabla,\;\;\; {\mathcal E}^1\circ{\mathcal E}.
\end{equation}
As 
 the curvature of $\Hat\nabla$ is given by
  ${\mathcal K}_{\Hat\nabla}=\Hat\nabla^1\circ\Hat\nabla$,
	we have the following proposition.
	
  \begin{Prop}\label{KHatNabla}
  The curvature of $\Hat\nabla=\nabla+{\mathcal E}$ is
  \begin{equation}\label{FormulaProp}
	{\mathcal K}_{\Hat\nabla}={\mathcal K}_{\nabla}+\nabla^1\circ{\mathcal E}+{\mathcal E}^1\circ{\nabla}  +{\mathcal E}^1\circ{\mathcal E}.
	\end{equation}
  \end{Prop}

The presheaf
 $$U\mapsto {\mathscr S}(U):= {\mathscr A}^1{}(U)\otimes_{\mathscr A{}(U)}{\rm Hom}_{\mathscr A{}|_U}({\mathscr F}|_U,\,{\mathscr F}|_U),$$
defines the ${\mathscr A}{}$-module ${\mathscr S}^+={\mathscr A}^1{}\otimes_{\mathscr A{}}{\mathscr Hom}_{\mathscr A{}}({\mathscr F},\,{\mathscr F})$  \cite[p. 231-232]{G-R}.
On the other hand, the ${\mathscr A}{}$-module ${\mathscr Hom}_{{\mathscr A}{}}\big({\mathscr F},\,{\mathscr A}{}^1\otimes_{{\mathscr A}{}}{\mathscr F}  \big)$ is the sheaf associated to the presheaf 
$$U\mapsto {\mathscr S}'(U):=
{\rm Hom}_{{\mathscr A}{}|_U}\big({\mathscr F}|_U,\,\big({\mathscr A}{}^1\otimes_{{\mathscr A}{}}{\mathscr F})|_U  \big).$$

Given $C=\alpha\otimes f\in {\mathscr S}(U)$, we define $C'\in {\mathscr S}'(U)$ as follows.  For any open $V\subset U$ and $\sigma\in {\mathscr F}(V)$,   we set 
\begin{equation}\label{aa'}
C'(\sigma)=\alpha|_V\otimes f|_V(\sigma).
\end{equation}
  Thus, the correspondence $C\mapsto C'$ is a morphism of presheaves 
	$\Phi: {\mathscr S}\to{\mathscr S}'$, which in turn determines
a map between the spaces of sections of the respective sheaves
\begin{equation}\label{mathcalQ}
 E\in\Gamma\big(M,\,   {\mathscr A}^1{}\otimes_{\mathcal A{}}{\mathscr Hom}_{\mathcal A{}}({\mathscr F},\,{\mathscr F})\big)\mapsto {\mathcal E}\in
{\rm Hom}_{{\mathscr A}{}}\big({\mathscr F},\,{\mathscr A}{}^1\otimes_{{\mathscr A}{}}{\mathscr F
 }  \big).
\end{equation}
(Just as the natural morphism $Z\otimes{\rm End}(P)\to{\rm Hom}(P,\,Z\otimes P)$ in the category of $R$-modules is not in general an isomorphism, (\ref{mathcalQ}) is not an isomorphism either).

 If the restriction of $E\in \Gamma\big(M,\,   {\mathscr A}^1\otimes_{\mathcal A{}}{\mathscr Hom}_{\mathcal A}({\mathscr F},\,{\mathscr F})\big)$  to an open subset $V$ of $M$
can be written as 
 \begin{equation}\label{E=eta*}
 { E}|_V=\sum_a\eta_a\otimes e_a, 
 \end{equation} 
 with $\eta_a\in {\mathscr A}^1{}(V)$ and
 $e_a\in {\rm Hom}_{{\mathscr A}|_V}({\mathscr F}|_V,\, {\mathscr F}|_V)$, then    
 \begin{equation}\label{BEsigma}
{\mathcal E}(\sigma)=\sum_a\eta_a\otimes e_a(\sigma),
\end{equation}
 for $\sigma\in{\mathscr F}(V)$.

We analyze the form of the operators involved in Proposition \ref{KHatNabla}
 when 
  ${\mathcal E}$ is the image of ${ E}$ by the map (\ref{mathcalQ}). Assumed that ${ E}|_V$ can be expressed as in (\ref{E=eta*}), 
  by (\ref{mathcalBk}), 
   $$  \big({\mathcal E}^1\circ{\mathcal E}\big)(\sigma)=\sum_{a,b}\eta_a\wedge\eta_b\otimes e_a(e_b(\sigma)).$$
   That is, on ${\mathscr F}(V)$,
	\begin{equation}\label{BEBE}
	{\mathcal E}^1\circ{\mathcal E}=\sum_{a,b}(\eta_a\wedge\eta_b)\otimes(e_a\circ e_b)=:E\wedge E.
	\end{equation}
	
   Similarly, by (\ref{nabla_k}),
   \begin{equation}\label{nabla1circ}
	\nabla^1\circ{\mathcal E}=\sum_a\big(d\eta_a\otimes e_a(\,.\,)-\eta_a\wedge\nabla(e_a(\,.\,))\big),
	\end{equation}
   on ${\mathscr F}(V)$. 
   
   If, for the sake of simplicity, we assume that $E=\eta\wedge e$ and that $\nabla (\sigma)$ can written as $\alpha\otimes \tau$, then ${\mathcal E}^1\circ \nabla(\sigma)=(\eta\wedge\alpha)\otimes\tau$. Thus, 
	\begin{equation}\label{E1nabla}
	{\mathcal E}^1\circ \nabla=E\wedge \nabla.
			\end{equation}
 
 \medskip
 \noindent
{\sc Connections on locally free sheaves.}
Let $F$ be a 
$C^{\infty}$ vector bundle of rank $m$  over $M$. By ${\mathscr F}$ we denote  the ${\mathscr A}{}$-module   of sections of $F$.  
 Let $\{u_i\,| \, i=1,\dots, n\}$ and $s=(s_1,\dots,s_m)$ be a  local frames for ${\mathscr A}{}^1$ and ${\mathscr F}$, respectively, defined  on  an open neighborhood $U$ of a point $x\in M$.
 Given an element $C'$ of ${\mathscr S}'(U)$, it can be written
on any $W\subset V$,
 as
$$C'=\sum_{a,i,b} f_{iab}  (u_i\otimes s_a)\otimes s^b,$$ 
where  $\check{s}=(s^1,\dots, s^m)$ is the frame dual of $s$,
and the  $f_{iab}$ are functions on $W$. Thus, $C'$ determines 
the element $C\in{\mathscr S}(W)$
$$C=\sum_{a,b} f_{ab}\otimes ( s_a\otimes s^b),$$
where $f_{ab}$ is the $1$-form $f_{ab}=\sum_if_{iab}u_i$. Hence,  $\Phi(W):{\mathscr S}(W)\to{\mathscr S}'(W)$ is an isomorphism. That is, $\Phi_{x}$ the map between the corresponding stalks  is also an isomorphism; i.e. ${\mathscr S}\simeq{\mathscr S}'$.  
Similarly,
\begin{equation}\label{isomo_sheaves}
  {\mathscr A}^k{}\otimes_{\mathcal A{}}{\mathscr Hom}_{\mathcal A{}}({\mathscr F},\,{\mathscr F})\big)
 \simeq
{\mathscr Hom}_{{\mathscr A}{}}\big({\mathscr F},\,{\mathscr A}{}^k\otimes_{{\mathscr A}{}}{\mathscr F
 }  \big).
\end{equation}
Thus, in particular, we have the following lemma. 

\begin{Lem}\label{L:1} 
If ${\mathscr F}$ is a locally free ${\mathscr A}$-module, then map  defined in (\ref{mathcalQ}) is an isomorphism of vector spaces.
\end{Lem}
Analogously, the curvature ${\mathcal K}_{\nabla}$ of a connection on ${\mathscr F}$  determines a vector
 of $\Gamma\big(M,\,{\mathscr A}{}^2\otimes_{{\mathscr A}{}}{\mathscr End}_{{\mathscr A}{}}({\mathscr F})\big)$, which will denoted $K_{\nabla}$. 

\smallskip

 As it is well-known, the connection  $\nabla$ induces a connection on   ${\mathscr End}_{\mathscr A}({\mathscr F})$, the sheaf of endomorphisms of ${\mathscr F}$,  by
\begin{equation} \label{nablah}
\nabla ({ h})=\nabla\circ{ h}- { h}\circ\nabla,
\end{equation}
for $h$ section of ${\mathscr End}_{\mathscr A}({\mathscr F}).$
 According to (\ref{nabla_k}) this connection defines in a natural way    an operator 
$$\nabla:{\mathscr A}{}^k\otimes_{{\mathscr A}{}}{\mathscr End}_{\mathscr A}({\mathscr F})
\to {\mathscr A}{}^{k+1}\otimes_{{\mathscr A}{}}{\mathscr End}_{\mathscr A}({\mathscr F}).$$
 On the $k$-form ${\mathscr End}_{\mathscr A}({\mathscr F})$-valued $\beta\otimes {h}$
\begin{equation}\label{nablabetah}
\nabla(\beta\otimes h)={d}\beta \otimes h +(-1)^k\beta\wedge\big( \nabla\circ h-h\circ\nabla  \big).
\end{equation}
  
If $E:=\beta\otimes h$ is a $1$-form, then it follows from (\ref{nabla1circ}) together with (\ref{E1nabla})  that 
\begin{equation}\label{nabla(E)}
\nabla(E)={\mathcal E}^1\circ \nabla+\nabla^1\circ {\mathcal E}.
\end{equation}

From Proposition \ref{KHatNabla} together with (\ref{BEBE}) and (\ref{nabla(E)}),
it follows the following proposition.
\begin{Prop}\label{P:2}
Let ${\mathscr F}$ be the sheaf of sections of a vector bundle $F$ over $M$. Given  and $E\in \Gamma\big(M,\, {\mathscr A}^1{}\otimes_{{\mathscr A}{}}{\mathscr End}_{{\mathscr A}{}}({\mathscr F})\big)$ and $\nabla$ a connection on  ${\mathscr F}$, then the curvature $\Hat{ K}$ of the connection $\Hat\nabla=\nabla+{\mathcal E}$ is
$$\Hat{ K}={K}+\nabla(E)+ E\wedge E,$$
where ${ K}$ is the curvature of $\nabla$.
\end{Prop} 
 
{\it Remark.}
The connection $\Hat\nabla=\nabla+{\mathcal E}$ on ${\mathscr F}$ defines the corresponding one on ${\mathscr End}_{\mathscr A}({\mathscr F})$ in accordance
with (\ref{nablah}). Thus, with the above notations, $\Hat\nabla(h)=\nabla+[{\mathcal E},\,h].$ Furthermore, one has the extension of $\Hat\nabla$ to a map on ${\mathscr A}^k{}\otimes_{{\mathscr A}{}}{\mathscr End}_{\mathscr A}({\mathscr F})$. This last operator 
acting on ${\mathscr A}{}^k\otimes_{{\mathscr A}{}}{\mathscr End}_{\mathscr A}({\mathscr F})$  is the sum of $\nabla$, and a
 morphism of sheaves
\begin{equation}\label{Emap}{\mathscr A}^k{}\otimes_{{\mathscr A}{}}{\mathscr End}_{\mathscr A}({\mathscr F})\to{\mathscr A}^{k+1}{}\otimes_{{\mathscr A}{}}{\mathscr End}_{\mathscr A}({\mathscr F}),\;\;
{ D}\mapsto {E}\wedge{ D}-(-1)^k{ D}\wedge{ E}. 
\end{equation}
In particular,  if ${ E}=\eta\otimes e,$ 
 then this map takes the form
\begin{equation}\label{betaob}
\beta\otimes h\in {\mathscr A}^k{}\otimes_{{\mathscr A}{}}{\mathscr End}_{\mathscr A}({\mathscr F})\mapsto 
(\eta\wedge\beta)\otimes[e,\,h]\in  {\mathscr A}^{k+1}{}\otimes_{{\mathscr A}{}}{\mathscr End}_{\mathscr A}({\mathscr F}).
\end{equation}

\smallskip

\noindent
{\sc Hermitian connections.} 
Let $F\to M$ be a $C^{\infty}$ a Hermitian vector bundle of rank $m$ over the manifold $M$.
We set ${\mathfrak u}({\mathscr F})$ for the sheaf of   antihermitian endomorphisms of ${\mathscr F}$. Obviously, ${\mathfrak u}({\mathscr F})$ is a ${\mathscr R}{}$-module, ${\mathscr R}{}$ being the sheaf of ${\mathbb R}$-valued $C^{\infty}$ functions on $M$.

The metric $\langle\,\cdot\,,\, \cdot\,\rangle$ on $F$ determines a metric on $\mathfrak{u}(F)$ which satisfies for all $x\in X$
\begin{equation}\label{langeadj}
\langle{\rm ad}_{c(x)}(a(x)),\,b(x)\rangle{} + \langle a(x),\,{\rm ad}_{c(x)}(b(x))\rangle{}=0,
\end{equation}
where  $a,\, b,\, c$ are any sections of  $\mathfrak{u}(F)$ and $a(x),\,b(x)$, and $c(x)$ are the corresponding vectors in the fibre
${\mathfrak u}({\mathcal F})(x)$ of the sheaf ${\mathfrak u}({\mathcal F})$ at $x$.

   A Hermitian connection $\nabla$ on ${\mathscr F}$ is a connection  compatible with the Hermitian structure; i.e. such that
\begin{equation}\label{Unitaryconnection}
\langle \nabla\sigma,\, \tau\rangle+ \langle \sigma,\, \nabla\tau\rangle=d\langle \sigma,\,\tau\rangle,
\end{equation}
for
 any sections $\sigma,\tau$ of ${\mathscr F}$. In this case, $\nabla$ determines a ${\mathbb R}_{X}$-morphism
 $$\nabla:{{\mathscr A}}^k{}\otimes_{{\mathscr R}{}}{\mathfrak u}({\mathscr F})\to {{\mathscr A}}^{k+1}\otimes_{{\mathscr R}{}}{\mathfrak u}({\mathscr F}).$$
 If furthermore  $\nabla$ is {\it integrable} (i. e. a flat gauge field),
 one has the following complex
$${\mathscr A}^{\bullet}{}\otimes_{{\mathscr R}{}}{\mathfrak u}({\mathscr F}):\;\;
{\mathfrak u}({\mathscr F})\overset{\nabla}{\longrightarrow} {\mathscr A}^{1}{}\otimes_{{\mathscr R}{}}{\mathfrak u}({\mathscr F})\overset{\nabla}{\longrightarrow}  {\mathscr A}^{2}{}\otimes_{{\mathscr R}{}}{\mathfrak u}({\mathscr F})  \to\dots$$
The Poincar\'e's lemma holds for this complex, and hence 
$${\mathscr H}^j\big({\mathscr A}^{\bullet}{}\otimes_{{\mathscr R}{}}{\mathfrak u}({\mathscr F})\big)=0,\;\;\; j>0.$$ 
 The $0$-cohomology is the local system ${\mathscr P}={\mathscr Ker}(\nabla)$, defined by the parallel sections of
 ${\mathfrak u}({\mathscr F})$.

The spectral sequence $E_2^{ij}=H^i\big(M,\,{\mathscr H}^j({\mathscr A}^{\bullet}{}\otimes_{{\mathscr R}{}}{\mathfrak u}({\mathscr F}))\big),$
converges to the hypercohomology ${\mathbb H}^*(X,\, 
{\mathscr A}^{\bullet}{}\otimes_{{\mathscr R}{}}{\mathfrak u}({\mathscr F}))$; that is, to the cohomology of the complex 
$${\mathscr A}^{\bullet}({\mathfrak u}(F)):=\Gamma(M,\, {\mathscr A}^{\bullet}{}\otimes_{{\mathscr R}{}}{\mathfrak u}({\mathscr F})).$$
 Hence,
 \begin{Prop}\label{P:3}
The cohomology groups of the sheaf ${\mathscr P}$ of parallel sections of ${\mathfrak u}({\mathscr F})$ satisfy
$$H^i(M,\,{\mathscr P})=H^i({\mathscr A}^{\bullet}({\mathfrak u}(F))).$$
\end{Prop}

\medskip

\subsection{Yang-Mills fields.}\label{Ss:Yang-Mills_fields}
 We assume now that $M$ is a Riemannian oriented manifold.
The   map between the sheaf of ${\mathscr A}^1{}$ of 
$C^{\infty}$ $1$-forms and the sheaf of ${\mathscr T}{}$ of $C^{\infty}$ vector fields induced by the metric will be denoted
$\lambda\mapsto \lambda^{\natural}.$   Given  $\xi,\zeta$ $k$-forms,   $(\xi,\,\eta)$ will denote the function defined by the relation
\begin{equation}\label{xiwedgestar}
\xi\wedge\star\zeta=(\xi,\,\zeta)\,{\rm dvol},
\end{equation}
where $\star$ is the Hodge star operator and ${\rm dvol}$ the volume form.

Given the vector field $v\in {\mathscr T}{}$, we set $\iota_v$ for the inner derivative defined by   $v$. One has the well-known   relation
\begin{equation}\label{iota}
(\lambda\wedge\xi,\,\zeta)=(\xi,\,\iota_{\lambda^{\natural}}(\zeta)).
\end{equation}

\smallskip

 As we said,  ${\mathscr A}^k({\mathfrak u}(F))$ will denote the vector space    
$${\mathscr A}^k({\mathfrak u}(F)):=\Gamma\big(M,\, {\mathscr A}^k{}\otimes_{{\mathscr R}{}} {\mathfrak u}({\mathscr F})\big)$$
 of global sections of 
 ${\mathscr A}^k{}\otimes_{{\mathscr R}{}} {\mathfrak u}({\mathscr F})$.

On the space ${\mathscr A}^k({\mathfrak u}(F))$ one defines the following product
\begin{equation}\label{langlelangle}
\langle\langle \xi\otimes\sigma,\,\xi'\otimes \sigma'\rangle\rangle:=\int_{M}(\xi\wedge\star\xi')\langle \sigma,\,\sigma'\rangle,
\end{equation}
 where $\langle\sigma,\,\sigma'\rangle$ is the function on $M$ defined by 
$\langle\sigma,\,\sigma'\rangle(x)=\langle \sigma(x),\,\sigma'(x)\rangle{}$ (here $\sigma(x)$ is the corresponding vector of 
the fibre of ${\mathfrak u}({\mathscr F})$).

Let $\nabla$  be an integrable  connection 
compatible with the metric of $F$. Then one has the complex 
\begin{equation}\label{elliptic-complex}
{\mathscr A}^0({\mathfrak u}(F))\overset{\nabla}{\longrightarrow}{\mathscr A}^1({\mathfrak u}(F))\overset{\nabla}{\longrightarrow}\dots 
\end{equation}
\begin{Lem}\label{L:elliptic}
The complex (\ref{elliptic-complex})
is elliptic.
\end{Lem}
{\it Proof.} Let $s$ be a unitary local frame of $F$ and $\check s$ denote its dual. An element $D\in {\mathscr A}^k({\mathfrak u}(F))$ can be written 
$$D=\sum_{i,j} {\sf D}_{ij}(s_i\otimes s^j)=s\cdot {\sf D}\cdot\check s,$$
where ${\sf D}$ is a matrix of $k$-forms. If ${\sf A}$ is the matrix of $\nabla$ in the frame $s$; that is,
$\nabla s=s\cdot{\sf A}$, and hence $\nabla\check s=-{\sf A}\cdot \check s$, then
\begin{equation}\label{nablaEnd}
\nabla(D)=s\cdot\big(d {\sf D}+{\sf A}\wedge {\sf D}-(-1)^k {\sf D}\wedge {\sf A} \big)\cdot\check s.
\end{equation}
That is, the operator $\nabla$ acting on the elements of   ${\mathscr A}^k({\mathfrak u}(F))$ 
is 
$$d+{\sf A}\wedge(\,\cdot\,)-(-1)^k(\,\cdot\,)\wedge{\sf A}.$$
Thus, the principal symbol of the operator $\nabla$ is equal to the one of the exterior differential operator $d$.
\qed

 Hence,
$$H^i({\mathscr A}^{\bullet}({\mathfrak u}(F)))={\rm Ker}(\Delta_i),$$
 where $\Delta_i$ is the corresponding Laplacian operator \cite[Chap. IV, Sect. 5]{Wells}.
From Proposition \ref{P:3},  it follows the   proposition.
\begin{Prop}\label{P:4} If  $\nabla$ is a connection flat, then 
 with the above notations
$$H^i (X,\,{\mathscr P})={\rm Ker}(\Delta_i).$$
\end{Prop}

\smallskip

Given a point $x\in M$, we denote by $\big({\mathscr A}^2{}\otimes_{{\mathscr R}{}}\mathfrak{u}({\mathscr F})\big)(x)$ the fibre of the sheaf
 ${\mathscr A}^2{}\otimes_{{\mathscr R}{}}\mathfrak{u}({\mathscr F})$ at $x$. 
 Let $\nabla$ be a connection ${\mathscr F}$ compatible with the metric $\langle\,,\,\rangle$. 
  $K_{\nabla}(x)$ will denote the vector of  $\big({\mathscr A}^2{}\otimes_{{\mathscr R}{}}\mathfrak{u}({\mathscr F})\big)(x)$ defined
	by ${ K}_{\nabla}$. 

The Yang-Mills functional is the map defined on the space of   Hermitian connections  by 
$$\nabla\mapsto {\sf YM}(\nabla)=  ||K_{\nabla}||^2,$$
where $||\;||$ is the norm determined by $\langle\langle\,,\,\rangle\rangle$ \cite[page 417]{Hamilton} \cite[page 44]{Moore} \cite[page 357]{Naber}. 
The Yang-Mills fields are the $\nabla$'s on which  ${\sf YM}$ takes  stationary values.

\smallskip

If $\nabla$ is a  connection on the Hermitian bundle $F$ compatible with the metric $\langle\,,\,\rangle$, from
 (\ref{Unitaryconnection}) together with (\ref{Leibniz}), it follows that
any other connection compatible with this metric
 has the form
$\Hat\nabla=\nabla+{\mathcal E}$, with   ${ E}\in{\mathscr A}^1{}(\mathfrak{u}({\mathscr F})).$

A ``variation'' of the connection $\nabla$ can be written  as $\nabla+\epsilon {\mathcal E}$.
By Proposition \ref{P:2},
 $K_{\epsilon}:=K_{\nabla+\epsilon {\mathcal E}}=K_{\nabla}+\epsilon\nabla(E)+\epsilon^2E\wedge E$. 
 Hence,  
$$\frac{d\,||K_{\epsilon}||^2}{d\epsilon}\Big|_{\epsilon=0}=2\langle\langle\nabla E,\, K_{\nabla}\rangle\rangle=2\langle\langle E,\, \delta K_{\nabla}\rangle\rangle,$$
where $\delta$ is the adjoint of   $\nabla$, when $\nabla$ is considered as an operator
\begin{equation}\label{nablaOperator}
\nabla:{\mathscr A}^1({\mathfrak u}({\mathscr F}))\to {\mathscr A}^2({\mathfrak u}({\mathscr F})). 
\end{equation}
Thus,
 $\nabla$ is a Yang-Mills field iff it satisfies the {\it Yang-Mills equation}
\begin{equation}\label{YMequation}\delta_{\nabla}(K_{\nabla})=0.
\end{equation}

  Given  ${E}\in{\mathscr A}^1({\mathfrak u}({\mathscr F}))$,  the 
	morphism of ${\mathscr A}{}$-modules (\ref{Emap}) gives rise to an operator  
   ${\mathscr A}^k({\mathfrak u}(F))\to {\mathscr A}^{k+1}({\mathfrak u}(F))$. 
 We will determine $\Hat\delta$, the formal adjoint operator of 
$$\Hat\nabla=\nabla+E:{\mathscr A}^1({\mathfrak u}({\mathscr F}))\to {\mathscr A}^2({\mathfrak u}({\mathscr F})).$$  
Obviously, $\Hat\delta=\delta+E^{\dagger}$, where 
 $E^{\dagger}$ is the adjont 
of the operator $E$.

If $E=\eta\otimes e\in{\mathscr A}^1({\mathfrak u}(F))$, by (\ref{betaob}), 
\begin{equation}\label{E:Omega1}
 E(\beta\otimes b)= (\eta\wedge \beta)\otimes [e,\,b],
\end{equation}
  for $\beta\otimes b\in {\mathscr A}^1({\mathfrak u}(F))$.
Given $\gamma\otimes c\in{\mathscr A}^2({\mathfrak u}(F))$, by (\ref{langlelangle}) 
$$\langle\langle E(\beta\otimes b),\,\gamma\otimes c\rangle\rangle=\int_{M}(\eta\wedge\beta\wedge\star\gamma)\langle[e,b],\,c\rangle.$$
From (\ref{xiwedgestar}) together with (\ref{iota}), it follows
$$(\eta\wedge\beta)\wedge\star\gamma= \beta\wedge\star(\iota_{\eta^{\natural}}(\gamma)).$$
On the other hand, by (\ref{langeadj}) 
$$\langle[e,\,b],\,c\rangle =\langle b,[c,\,e]\rangle.$$
Thus,
$$\langle\langle E(\beta\otimes b),\,\gamma\otimes c\rangle\rangle=
\int_M\beta\wedge\star(\iota_{\eta^{\natural}}(\gamma))\langle b,[c,\,e]\rangle=\langle\langle \beta\otimes b,\,E^{\dagger}(\gamma\otimes c\rangle\rangle,$$
where
$$E^{\dagger} (\gamma\otimes c)=\iota_{\eta^{\natural}}(\gamma)\otimes [c,\,e].$$
We have proved the following proposition.
\begin{Prop}\label{P:Edagger}
Given $E=\eta\otimes e\in{\mathscr A}^1({\mathfrak u}(F))$ the adjoint of the operator $E$ defined in (\ref{E:Omega1}) is
$$E^{\dagger}:{\mathscr A}^2({\mathfrak u}(F))\to {\mathscr A}^1({\mathfrak u}(F)), \;\;\;\;\gamma\otimes c\mapsto \iota_{\eta^{\natural}}(\gamma)\otimes [c,\,e].$$
\end{Prop}

As a consequence of (\ref{YMequation}) and Proposition \ref{P:2}, we have the following theorem.
\begin{Thm}\label{Th:YMHat}
Let $\nabla$ be a integrable connection on $F$ compatible with the Hermitian structure and 
${ E}\in{\mathscr A}^1({\mathfrak u}({\mathscr F}))$. Then $\Hat\nabla=\nabla+{\mathcal E}$ is a Yang-Mills connection iff 
$$\delta(\nabla(E))+\delta(E\wedge E)+E^{\dagger}(\nabla(E))+ E^{\dagger}(E\wedge E)=0,$$
where $E^{\dagger}$ is the operator defined in Proposition \ref{P:Edagger}, and $\delta$ is the adjoint operator  (\ref{nablaOperator}).
\end{Thm}

If $F$ admits integrable connections, on these gauge fields, ${\sf YM}$ takes the minimal value. The vacuum states of corresponding Yang-Mills theory are the gauge equivalence classes integrable connections.  The space of vacuum states is in bijective correspondence to the set of equivalence classes of the representations of $\pi_1(M)$ in ${\rm U}(m)$ \cite{Vina20}.

Not every vector bundle $F$ on the manifold $M$ admits an integrable connection.
By the Frobenius theorem, the existence of an integrable connection  on ${F}$ is equivalent to the fact that   $F$ admits a family of local frames, whose domains cover $M$ and such that the corresponding transition functions are constant; i.e., $F$ is a flat vector bundle 
\cite[p. 5]{Koba}.

	\smallskip
	
Assumed that $F$ is a flat vector bundle,  in the next section, we will study properties of     continuous families $\{\nabla_t\}_t$ of Yang-Mills fields with initial point at a vacuum state. That is, properties of the evolution of a Yang-Mills field from a vacuum state.


\subsection{Curves of Yang-Mills connections}\label{Ss:curves_Y-M}

Throughout  this section $\nabla$ will denote a {\it flat } Hermitian connection on $F$. 
 A curve of Hermitian connections on $F$ with initial point at $\nabla$ has the form
\begin{equation}\label{nablat}
\nabla_t=\nabla+{\mathcal E}(t),\;\;\;\ t\geq 0,
\end{equation}
with ${ E}(t)\in{\mathscr A}^1({\mathfrak u}({\mathscr F}))$ and ${ E}(0)=0$.
We define the following vectors of ${\mathscr A}^1({\mathfrak u}({\mathscr F}))$ 
$$E_1:=\lim_{t\to 0}E(t)/t,\;\;\; E_2:=\lim_{t\to 0}(E(t)-tE_1)/t^2.$$
Thus, $E(t)=tE_1+t^2E_2+O(t^3)$. We set $C_E:=\nabla E_2+ E_1\wedge E_1 \in {\mathscr A}^2({\mathfrak u}({\mathscr F}))$.
 \begin{Prop}\label{varphiE}
If $\nabla_t$ is a curve of Yang-Mills connections, then
\begin{enumerate} 
\item $E_1$ determines an element $[E_1]\in H^1(X,\,{\mathscr P})$,
 \item   $C_E$ defines  an element of $H^2(X,\, {\mathscr P})$.
 \end{enumerate}
\end{Prop}

{\it Proof.}
  By Theorem \ref{Th:YMHat},
$$0=t\delta(\nabla E_1)+t^2\delta(\nabla E_2)+t^2\delta(E_1\wedge E_1)+t^2E^{\dagger}_1(\nabla E_1)+O(t^3).$$
Thus,
\begin{equation}\label{deltanablaE1}
\delta(\nabla E_1)=0,\;\;\delta(\nabla E_2)+\delta(E_1\wedge E_1)+E^{\dagger}_1(\nabla E_1)=0.
\end{equation}

The condition $\delta(\nabla E_1)=0$ implies
$$0=\langle\langle \delta(\nabla E_1),\,E_1\rangle\rangle=\langle\langle  \nabla E_1,\,\nabla E_1\rangle\rangle; $$
that is $\nabla(E_1)=0$. From Proposition \ref{P:3}, it follows the first assertion of the proposition.

Consequently,  the second equation in (\ref{deltanablaE1}) reduces to
$$\delta\big(\nabla(E_2)+E_1\wedge E_1\big)=0.$$
That is, $\delta(C_E)=0$. Furthermore, as $\nabla$ is flat,
$$\nabla(C_E)=\nabla^2(E_2)+\nabla(E_1)\wedge E_1-E_1\wedge\nabla(E_1)=0.$$
From $\delta(C_E)=0$ and $\nabla(C_E)=0$, together with Proposition \ref{P:4}, it follows that $C_E$ defines an
 element of $H^2(X,\,{\mathscr P})$. \qed 

\begin{Prop}\label{CEflat}
If $\nabla_t=\nabla+{\mathcal E}(t)$   is a curve of integrable connections, then 
 $C_E=0$.
\end{Prop}
{\it Proof.} As $\nabla_t$ is flat, by Proposition \ref{P:2},
$$\nabla E(t)+E(t)\wedge E(t)=0.$$ Hence
$$t\nabla E_1+t^2\big(\nabla E_2+E_1\wedge E_1\big)+O(t^3)=0.$$
Thus $\nabla E_1=0$ and $C_E=0$.        \qed

\smallskip

\noindent
{\it Remark.}
By ${\mathscr U}({\mathscr F})$ we denote the subsheaf of ${\mathscr Aut}({\mathscr F})$ consisting of the unitary automorphisms. Those automorphisms are the gauge transformations that preserve the Hermitian metric. Given ${G }\in\Gamma( X,\,{\mathscr U}({\mathscr F}))$, the composition 
${G}\circ\nabla\circ{G}^{-1}$ defines a new connection gauge equivalent to $\nabla$.  A smooth family ${G}_t$ of
 elements of $\Gamma(X,\,{\mathscr U}({\mathscr F}))$, where ${G}_0$ is the identity, gives rise to a development
$$G_t={\rm id}+tA_1+t^2A_2+O(t^3).$$


Taking into account that ${ G}_t^{-1}={\rm id}-tA_1+t^2(A_1A_1-A_2)+O(t^3),$ it is straightforward to check that the corresponding $E_1$ and $E_2$ associated to the curve of connections ${G}_t\circ\nabla\circ{G}_t^{-1}$ are 
$$E_1=-\nabla(A_1),\;\; E_2=\nabla(A_1A_1)-\nabla(A_2)-A_1\nabla(A_1).$$
Therefore for this curve $[E_1]=0$ and of course $C_E=0$.
\smallskip


\subsection{Curve connecting vacuum states}\label{Ss:Connecting} 
Let $T=T^2$ be the $2$-torus $\big[0,\,1]\times[0,\,1]\big)/\sim,$ with
$$(x,\,0)\sim (x,\,1),\;\;(0,\,y)\sim (1,\,y),  $$
endowed with the metric $dx^2+dy^2$. We denote by $w=dx\wedge dy$ the corresponding volume form.  Let $F=T\times {\mathbb C}^2\to T$ be the trivial bundle equipped with the obvious Hermitian structure. On $F$ we consider the connection $\nabla$ defined by the exterior derivative.   That is, $\nabla s_i=0$, $i=1,2$, where $s_1(z)=(z,\,(1,0))$ and $s_2(z)=(z,\,(0,1))$, for $z\in T^2$. Evidently, $\nabla$  is a connection compatible with the Hermitian structure.

We will construct  a new connection on $F$, $\Hat\nabla=\nabla+ {\mathcal E}$, where
 $E\in {\mathscr A}^1_{T} \otimes_{{\mathscr A}_T}\mathfrak{su}(F)$ is defined as follows.
By $\sigma_a$, $a=1,2,3$ we denote the Pauli matrices. Then $e_a(z)=(z,\, i\sigma_a)$,  is a section of $T\times\mathfrak{su}(2)\to T$. We set  
\begin{equation}\label{E=alpha}
E:=\alpha\otimes e_1+\beta \otimes e_2+\gamma\otimes  e_3,
\end{equation}
with $\alpha,\beta,\gamma$ $1$-forms on $T$. As $[e_a,\,e_b]=-2\epsilon_{abc}e_c$, 
then (see (\ref{BEBE}))
\begin{equation}\label{EwedgeE}
E\wedge E=-2\big( (\beta\wedge\gamma)\otimes e_1 + (\gamma\wedge\alpha)\otimes e_2 + (\alpha\wedge\beta)\otimes e_3 \big).
\end{equation}
We will assume that $E\wedge E=0$; that is, $\beta\wedge\gamma = \gamma\wedge\alpha = \alpha\wedge\beta=0$.
Since $d\alpha,\,d\beta,\,d\gamma$ $2$-form, there exist functions $h_1,h_2,h_3$ such that
\begin{equation}\label{dalph=0h}
d\alpha=h_1w,\; d\beta=h_2w,\;d\gamma=h_3w.
 \end{equation}
 Hence,
\begin{equation}\label{nablaE=}
\nabla E=w\otimes\sum_{a=1}^3h_ae_a.
\end{equation}

By Theorem \ref{Th:YMHat},  $\Hat\nabla=\nabla+{\mathcal E}$ is a Yang-Mills connection iff
\begin{equation}\label{deltanabla E}
\delta(\nabla E)+E^{\dagger}(\nabla E)=0.
 \end{equation}
From (\ref{nablaE=}), one obtains   
\begin{equation}\label{deltanablaE}
\delta(\nabla E)=-\star d\star(\nabla E)=-\sum_a \star dh_a\otimes  e_a.
\end{equation}
From Proposition \ref{P:Edagger} together with (\ref{nablaE=}), it follows
$$(\alpha\otimes e_1)^{\dagger}(\nabla E)=\sum_a\iota_{\alpha^\sharp}(h_aw)\otimes[e_a,\,e_1]=
2\big(  \iota_{\alpha^\sharp}(h_2w)\otimes e_3- \iota_{\alpha^\sharp}(h_3w)\otimes e_2 \big).$$
Expressions for $(\beta\otimes e_2)^{\dagger}(\nabla E)$ and for $(\gamma\otimes e_3)^{\dagger}(\nabla E)$ can be obtained in a similar way.
From (\ref{E=alpha}), we deduce  
 \begin{align}
\notag &E^{\dagger}(\nabla E)=2\Big(
  \big(-\iota_{\gamma^\sharp}(h_2w) + \iota_{\beta^\sharp}(h_3w)\big)\otimes e_1 +\\ \notag
	&+\big( -\iota_{\alpha^\sharp}(h_3w) + \iota_{\gamma^\sharp}(h_1w)\big)\otimes e_2  
	+ \big(-\iota_{\beta^\sharp}(h_1w)+ \iota_{\alpha^\sharp}(h_2w)\big)\otimes e_3\Big). \notag 
	\end{align} 
	Hence, $\Hat\nabla$ is a Yang-Mills connection if
	\begin{align}
	 \star d h_1=2\big(\iota_{\gamma^\sharp}(h_2w)- &\iota_{\beta^\sharp}(h_3w)   \big),\;\;
	\star dh_2=2\big(  \iota_{\alpha^\sharp}(h_3w) - \iota_{\gamma^\sharp}(h_1w)\big)      \big) \label{align12}\\
&\star dh_3=2\big( (\iota_{\beta^\sharp}(h_1w)- \iota_{\alpha^\sharp}(h_2w)  \big). \label{align3}
\end{align}

If $\alpha=a_1dx+a_2dy,\; \beta=b_1dx+b_2dy,\; \gamma=c_1dx+c_2dy$, then
$$\iota_{\alpha^{\sharp}}(w)=a_1dy-a_2dx,\; \iota_{\beta^{\sharp}}(w)=b_1dy-b_2dx,\;\iota_{\gamma^{\sharp}}(w)=c_1dy-c_2dx.$$
When $\gamma=0$, but $\alpha\ne 0\ne\beta$ then $h_3=0$ and from (\ref{align12}), 
 we deduce that $h_1$ and $h_2$ are constant functions, with $h_1\ne 0\ne h_2$. From (\ref{align3}), it follows
$$\frac{h_2}{h_1} =\frac{b_1}{a_1}=\frac{b_2}{a_2};$$
that is, $\beta=\lambda\alpha$, where $\lambda$ is a non zero constant. Thus, $E=\alpha\otimes (e_1+\lambda e_2)$ defines a Yang-Mills field $\Hat\nabla$.
From (\ref{nablaE=}), it follows $ \nabla E\ne 0$.    From Proposition \ref{P:2}, one deduces the following proposition.
\begin{Prop}\label{P:Torus}
  With the above notations, the   connection $\Hat\nabla=\nabla+{\mathcal E}$, where $E=\alpha(e_1+\lambda e_2)$ with $\alpha\ne 0$ and $\lambda\ne 0$, is a non flat Yang-Mills connection. 
\end{Prop}

\smallskip

{\it Remark.} If $\beta=\gamma=0,$ the equations for the $h_a$ imply that $h_1$ is constant and $h_2=h_3=0$, obviously. From (\ref{dalph=0h}), it follows $h_1=0$ (otherwise $w$ would be exact), and by  (\ref{nablaE=}) $\nabla E=0$. That is, the connection $\Hat \nabla$ defined by $E=\alpha\otimes e_1$ is flat.
Moreover, if $\alpha$ is not exact, then $\nabla$ and $\Hat\nabla$ are not gauge equivalent (in fact, $\Hat\nabla$ defines a non trivial representation of   $\pi_1(T^2)$, but the representation associated to    $\nabla$ is the trivial one).

\smallskip
 For instance, $\tilde\alpha =(\sin\pi x )dy+(\cos\pi y)dx$ is a non closed  $1$-form on $T$ and the corresponding connection $\Hat\nabla$ defined by $E=\tilde\alpha\otimes (e_1+\lambda e_2)$  (with $\lambda\ne 0$) is not flat, by Proposition \ref{P:Torus}. 
For $t\in[0,\,1]$, we define 
\begin{equation}\label{Et=t}
E_t=t\big((\sin\pi x )dy+(\cos\pi y)dx\big)\otimes(e_1+(1-t)e_2).
\end{equation}
We have a curve $\Hat\nabla_t=\nabla+{\mathcal  E}_t$ of connections, where $\Hat\nabla_0$ is obviously flat.  $E_1$ reduces to $\tilde\alpha\otimes e_1$; thus, according to the Remark, $\Hat\nabla_1$ is a flat gauge field. Hence, one has a curve on Yang-Mills fields 
 connecting two different vacuum states of the Yang-Mills theory. From Proposition \ref{P:Torus}, it follows: 
\begin{Prop}\label{curveconn}
The family $E_t$ defined in (\ref{Et=t})  determines a curve $\{\Hat\nabla_t\}_{t\in[0,\,1]}$ of Yang-Mills fields  on $F=T\times {\mathbb C}^2$ satisfying
\begin{enumerate}
\item $\Hat\nabla_0$ and $\Hat\nabla_1$ are gauge inequivalent integrable connections 
\item For each $t\in(0,\,1)$, $\Hat\nabla_t$ is a non flat   Yang-Mills connection. 
\end{enumerate} 

\end{Prop}
 

\section{Yang-Mills fields on  coherent sheaves}\label{S:CoherentSheaves}
\subsection{Holomorphic gauge fields on vector bundles}\label{Jet1}
 In this  section $X$ is a complex analytic manifold, ${\mathscr O}{}$ will denote the sheaf of holomorphic functions on $X$, and $\Omega^k{}$ the sheaf of holomorphic $k$-forms.
Let $\tilde\pi:F\to X$ be a {\it holomorphic}  vector bundle over   $X$. By $F$ will be also denoted the   sheaf of their sections. 
We review the construction of the    $1$-jet bundle $J^1(F)$   to define the  holomorphic gauge fields on $F$ and to show the obstruction to the existence of these fields.

  We denote by $\Gamma(q)$   the set of sections of $F$ defined in some neighborhood of the point $q\in X$. In $\Gamma(q)$ one defines the following equivalence relation
$\sigma\sim\tau$, iff 
$$\sigma(q)=\tau(q)\;\;\hbox{and}\;\; \sigma_*|_{T_qX}=\tau_*|_{T_qX},$$
where 
$$\sigma_*|_{T_qX}:T_qX\longrightarrow T_{\sigma(q)}F,$$
is the linear map induced between the tangent spaces.

 The equivalence class of $\sigma$ is denoted $j^1_q(\sigma)$. One has $J^1(F)$, the first jet bundle of $F$,
$$J^1(F)=\{j^1_q(\sigma)\,|\, q\in X,\;\sigma\in\Gamma(q)\}\overset{p}{\longrightarrow} X,\;\;\; j^1_q(\sigma)\mapsto q.$$ 

Given  $\sigma$  a section of $F$ defined on $U$, it  determines a section $j^1(\sigma)$ of $J^1(F)\to X$ defined by $j^1(\sigma)(q)=j^1_q(\sigma)$. If $f$ is a holomorphic function on $U$, one defines
$f\cdot j^1(\sigma):=j^1(\tau)$, where $\tau$ is a section of $F$, such that, $\tau(q)=\sigma(q)$ and $\tau_*|_{T_qX}=f(q)\sigma_*|_{T_qX}$ for all $q\in U$. 
If $\tau'$  satisfies also the preceding conditions, then $j^1(\tau)=j^1(\tau')$; thus, $f\cdot j^1(\sigma)$ is well defined.
 It is easy to check $(gf)\cdot j^1(\sigma)=g\cdot(f\cdot j^1(\sigma))$. 

 On the other hand, we  have the epimorphism of ${\mathscr O}$-modules. 
$$\pi:J^1(F)\to F,\;\;j^1_q(\sigma)\mapsto \sigma(q),$$
 and the corresponding short exact sequence of ${\mathscr O}$-modules
$$0\to {\rm ker}(\pi)\rightarrow J^1(F) \overset{\pi}{\rightarrow} F\to 0.$$

If $\pi(j^1_q(\sigma))=0$, then $\sigma_*|_{T_qX}$ maps $T_qX$ into the tangent vector space to the fibre $\tilde\pi^{-1}(q)$ at the point $0$.
 This last space can be identified with the fibre $F_q$. Thus,
 $\sigma_*|_{T_qX}$ defines an element of
$T^*_qX\otimes F_q=\Omega^1_q(F)$, the fibre at $q$ of the sheaf of $F$-valued holomorphic $1$-forms. Thus, one has
  the following short exact sequence
\begin{equation}\label{shortexact}
0\to\Omega^1(F)\overset{i}{\to} J^1(F)\overset{\pi}{\to} F\to 0.
\end{equation}


The natural map $j^1:F\to J^1(F)$ is a morphism of abelian sheaves, but not an arrow of ${\mathscr O}$-modules.
In fact,
\begin{equation}\label{jfsigma}
j^1(f\sigma)=i(df\otimes\sigma\big)+f\cdot j^1(\sigma).
\end{equation}
As $\pi\circ j^1={\rm id}$, (\ref{shortexact}) splits in the category of abelian sheaves
\begin{equation}\label{splitingJ1}
J^1(F)=F\oplus\Omega^1(F)\;\;\;\; \text{(as abelian sheaves).}
\end{equation}
 Expressing an element of $J^1(F)$ as $\sigma\oplus\alpha$, the product $f\cdot(\sigma\oplus\alpha)$ corresponds to
$$j^1(f\sigma)+i(f\alpha)=i(df\otimes \sigma)+f\cdot j^1(\sigma)+i(f\alpha),$$
 hence, by the identification (\ref{splitingJ1}), one can write
\begin{equation}\label{fsigma+alpha}
f\cdot(\sigma\oplus\alpha)=f\sigma\oplus\big(df\otimes\sigma+f\alpha  \big).
\end{equation}

On the other hand,   
if $D$ is a right inverse of $\pi$ in the category of ${\mathscr O}{}$-modules,  $\pi\circ D={\rm id}$, then 
$${\rm im}(j^1-D)\subset{\rm ker}(\pi)=i\big(\Omega^1{}(F)\big).$$
  For $\nabla:=j^1-D$ one has  by (\ref{jfsigma})
$$\nabla(f\sigma)=(j^1-D)(f\sigma)=i\big(df\otimes\sigma\big)+fj^1(\sigma)-fD(\sigma)=i\big(df\otimes\sigma\big)+f\nabla\sigma.$$

Thus, if $D$ is a right inverse of $\pi$, it defines a morphism of abelian sheaves
$$\nabla:F\to\Omega^1{}(F)$$
satisfying the Leibniz's rule for the product of holomorphic functions by sections. That is, $\nabla$ is a {\it holomorphic} connection. 

Hence, if there exists a holomorphic connection on $F$, then the exact sequence (\ref{shortexact}) splits, and conversely. In other words, the existence of a holomorphic connection on $F$ is equivalent to the vanishing in $ {\rm Ext}^1(F,\,\Omega^1(F))$ of  the extension defined by 
(\ref{shortexact}). The element of ${\rm Ext}^1(F,\,\Omega^1(F))$ determined by that exact sequence is called the Atiyah class of $F$, and it is denoted by $a(F)$. In particular, if $F$ is a line bundle, then  $a(F)=0$ whenever the first Chern class $c_1(F)$ vanishes \cite[Prop. 12]{Atiyah}

Since $F$ is a locally free ${\mathscr O}{}$-module, for any  coherent  ${\mathscr O}{}$-module ${\mathscr G}$ the sheaves 
${\mathscr Ext}^q(F,\,{\mathscr G})$ vanish for all $q>0$ and ${\mathscr Ext}^0(F,\,{\mathscr G})={\mathscr Hom} (F,\,{\mathscr G})$ \cite[Chap. III, Prop. 6.5]{Hartshorne}. On the other hand, the local to global  spectral sequence $E^{p,q}=H^p(X,\, {\mathscr Ext}^q( \cdot\,,\,\cdot))$ abuts to 
${\rm Ext}^{p+q}( \cdot\,,\,\cdot).$ Hence, ${\rm Ext}^{p}(F,\,{\mathscr G})=H^p(X,\,{\mathscr Hom}(F,\,{\mathscr G}))$. In particular, 
for ${\mathscr G}=\Omega^1{}\otimes_{{\mathscr O}{}} F$, one has
 $${\rm Ext}^{1}(F,\, \Omega^1{}\otimes_{{\mathscr O}{}} F )=H^1(X,\,\Omega^1{}\otimes_{{\mathscr O}{}} {\mathscr End}(F)).$$
 That is,  the Atiyah class $a(F)$ is an element  of the first cohomology of $\Omega^1{}\otimes {\mathscr End}(F)$.   By Cartan's Theorem B \cite[p. 243]{G-Ro}, we have the following proposition.  
\begin{Prop}\label{Stein} 
If $X$ is a Stein manifold, then any holomorphic vector bundle over $X$ admits holomorphic gauge fields.
\end{Prop} 

Let $F$ be a holomorphic vector bundle such that $a(F)=0$. If $X$ is a compact K\"ahler manifold and $F$, as $C^{\infty}$ bundle, is  Hermitian,   using the  $\star$ operator, one can define the corresponding Yang-Mills functional on the set holomorphic gauge fields compatible with the Hermitian structure. Let $\nabla$ be such a field. Any other field in this set is of the form $\nabla+E$, with 
$$E\in\Gamma(X,\,\Omega^1\otimes_{\mathscr O}\mathfrak{hu }(F)),$$ 
where $\mathfrak{hu}(F)$ is the sheaf of holomorphic antihermitian endomorphisms of $F$. Then  the Yang-Mills functional gives rise to the map
$$E\in\Gamma(X,\,\Omega^1\otimes_{\mathscr O}\mathfrak{hu }(F))\mapsto ||K_{\nabla+E}||^2.$$
The curvature $K_{\nabla+\epsilon E}=K_{\nabla}+\epsilon\nabla(E)+O(\epsilon ^2)$. Thus
$$\frac{d}{d\epsilon}\Big|_{\epsilon=0}||K_{\nabla+\epsilon E}||^2=2\langle\langle K_{\nabla},\,\nabla(E)\rangle\rangle.$$
Hence, $\nabla$ is a holomorphic Yang-Mills field iff $\langle\langle K_{\nabla},\,\nabla(E)\rangle\rangle=0$, for all $E$.

As  $\Omega^k\otimes_{\mathscr O}\mathfrak{hu }(F))$ is a coherent sheaf \cite[Annex. \S 4]{G-R}, $H^0(X,\, \Omega^k\otimes_{\mathscr O}\mathfrak{hu }(F))$ is a finite dimensional vector space \cite[page 700]{G-H}.    We denote by  $r_k:={\rm dim}\,H^0(X,\, \Omega^k\otimes_{\mathscr O}\mathfrak{hu }(F))$.

Let $\beta_1,\dots,\beta_{r_1}$ be a basis of $H^0(X,\, \Omega^1\otimes_{\mathscr O}\mathfrak{hu }(F))$ . Then the condition of being Yang-Mills reduces to the following conditions for the corresponding curvature:
 \begin{equation}\label{Reallange}
\langle\langle K_{\nabla},\,\nabla(\beta_i)\rangle\rangle=0,\;\; i=1,\dots,r_{1}.
\end{equation}

On the other hand, $K_{\nabla}$ is a vector of the finite dimensional space $H^0(X,\, \Omega^2\otimes_{\mathscr O}\mathfrak{hu }(F))$, which, according to (\ref{Reallange}), must be orthogonal to the image of the map 
\begin{equation}\label{map_nabla}
 \nabla: H^0(X,\, \Omega^1\otimes_{\mathscr O}\mathfrak{hu }(F))\to H^0(X,\, \Omega^2\otimes_{\mathscr O}\mathfrak{hu }(F)).
\end{equation}
 When $r_2\leq r_1$, we say that  $\nabla$ is {\it generic} if the map (\ref{map_nabla}) is surjective.
  Consequently, the curvature of any generic Yang-Mills connection must be zero.  One has the following proposition.
	
\begin{Prop}\label{P:generic} Let $F$ be a holomorphic Hermitian vector bundle on the compact K\"ahler manifold $X$. 
 If $r_2\leq r_1$, then all the generic holomorphic Yang-Mills fields   on $F$ are flat.
\end{Prop} 

When $F$ is a line bundle $r_2={\dim}\,H^0(X,\Omega^2)$ and $r_1={\dim}\,H^0(X,\Omega^1)$ are the Hodge numbers $h^{2,0}(X)$ and $h^{1,0}(X)$, respectively.
\begin{Cor}\label{Hodge}
 When $F$ is a line bundle, if $h^{2,0}(X)\leq h^{1,0}(X)$, then any generic holomorphic Yang-Mills connection on $F$ is flat.
\end{Cor}

Let $L$ be a Hermitian {\it line} bundle over the compact K\"ahler manifold $X$. Let us assume that $L$ is flat, thus   the space of holomorphic connections on $L$ is nonempty. 
Since 
$H^0(X,\, L)$ is a finite dimensional vector space, let $s_1,\dots,s_k$ be a basis of this space. After a normalization, we can assume that the basis is orthonormal, with respect to the product of sections defined by the Hermitian structure of $L$ and the K\"ahler structure of $X$. 

We denote by $\nabla_0$ the holomorphic connection determined by $\nabla_0(s_i)=0$ for $i=1,\dots,k$.
 If $\nabla$ is other holomorphic connection  on $L$, then 
$\nabla=\nabla_0+{ E}$, with 
$${ E}\in {\rm Hom}_{{\mathcal O}{}}\big(L,\,\Omega^1{}\otimes_{{\mathcal O}{}} L   \big)=\Gamma(X,\,\Omega^1{}).$$
The curvature of   $\nabla$ is ${ K}=\nabla_0({ E})$. 

The form ${ E}$ can be written as a linear combination 
$E=\sum_i\alpha_is_i$, with $\alpha_i$ holomorphic $1$-form. Hence, $K=\sum_i d(\alpha_i)s_i$. Setting  
 $$G(d\alpha_1,\dots,d\alpha_k):=||K||^2=\sum_i||d\alpha_i||^2,$$
the only stationary point of $G$ is $d\alpha_1=\dots=d\alpha_k=0$, which is defined by the connections whose curvature vanishes.
That is, we have the following proposition.

\begin{Prop}\label{noYangMills}
If $L$ is a Hermitian flat  line bundle on the compact K\"ahler manifold $X$,  then all the  holomorphic Yang-Mills fields on $L$ are flat.
\end{Prop}


\subsection{Holomorphic connections on coherent sheaves}\label{connections_on:Coherent}
There is another construction of the jet bundle which admits a translation to the case of coherent ${\mathscr O}{}$-modules. We summarize the development of \cite[p. 6]{Deligne}. By ${\mathscr I}$ we denote the ideal sheaf of the diagonal embedding $X\hookrightarrow  X\times X$. The structure sheaf of the first infinitesimal neighborhood of the diagonal is ${\mathcal O}_{X\times X}/ {\mathscr I}^2$ 
\cite[p. 5]{Grothendieck}. We set $X^{(1)}$ for the corresponding subscheme. Denoting by $p_1,p_2$ the projections $X\times X \rightrightarrows X$, one can define on ${\mathscr O}_{X^{(1)}} $ left and right ${\mathscr O}{}$-module structures via $p_1$ and $p_2$ respectively. 

One can think of sections of ${\mathscr O}{}$ as functions $h(x)$ and the sections of  ${\mathscr O}_{X^{(1)}} $ as classes of functions
 $g(x,\,y)$ modulo ${\mathscr I}^2$. The right and left ${\mathscr O}{}$-structues are defined by
$$((g +{\mathscr I}^2) \cdot f)(x,y)=g(x,y)f(y)+{\mathscr I}^2,\;(f\cdot (g+{\mathscr I}^2))(x,y)=f(x)g(x,y)+{\mathscr I}^2.$$

We have the   ${\mathbb C}$-linear map
$$m: {\mathscr O}_{X^{(1)}}\rightarrow{\mathscr O}{}\oplus \Omega^1{},\;\; g+{\mathscr I}^2\mapsto g(x,x)\oplus (d{}g)(x,x),$$
where $d{}g$ is the exterior differential of $g$ regarded as a function of $x$. It is easy to check that $m$ is a well defined morphism of abelian sheaves. Moreover, 
$$m(f\cdot (g+{\mathscr I}^2))=fm(g)+ \tilde g df,\;\;\; m((g+{\mathscr I}^2)\cdot f)=m(g) f$$
where $\tilde g(x)=g(x,\,x)$.
If we define a left ${\mathscr O}{}$-module structure on ${\mathscr O}{}\oplus \Omega^1{}$
by 
$$f\cdot(h\oplus\alpha):=fh\oplus(hdf+f\alpha),$$
then $m$ is an isomorphism of ${\mathscr O}{}$-modules.

In this way, for the locally free sheaf $F$, the ${\mathscr O}{}$-module ${\mathscr O}_{X^{(1)}}\otimes_{{\mathscr O}{}}F=F\oplus\Omega^1{}(F)$ with ${\mathscr O}{}$-module structure
$$f\cdot(\sigma\oplus \beta)=f\sigma\oplus(df\sigma+f\beta)$$
is isomorphic to the ${\mathscr O}{}$-module $J^1(F)$ (see (\ref{splitingJ1}) and (\ref{fsigma+alpha})).

\smallskip

It is possible to define the sheaf of $1$-jets of a coherent sheaf on $X$, generalizing the definition given for  vector bundles.
Let ${\mathscr F} $ be  a coherent ${\mathscr O}{}$-module,   one defines
\cite[p.6]{Deligne}
\begin{equation}\label{Def:J1F}
{\mathscr J}^1({\mathscr F}) :={\mathscr O}_{X^{(1)}}\otimes_{{\mathscr O}{}} {\mathscr F}.
\end{equation}
One has the short exact sequence of ${\mathscr O}{}$-modules
\begin{equation}\label{exatseq} 0\to\Omega^1{}\otimes_{{\mathscr O}{}}{\mathscr F}\overset{i}{\rightarrow} {\mathscr J}^1({\mathscr F}) \overset{\pi}{\rightarrow} {\mathscr F}\to 0.
\end{equation}

\begin{Def}  A holomorphic connection on  the coherent sheaf ${\mathscr F}$ is a right inverse of the morphism.
$\pi:{\mathscr J}^1({\mathscr F}) \to {\mathscr F}$
\end{Def}
That is, a connection on ${\mathscr F}$  is a splitting of (\ref{exatseq}). Then, as in Subsection \ref{Jet1}, 
there is a morphism of ${\mathbb C}_{X}$-modules   
$$\nabla:{\mathscr F}\to \Omega^1{}\otimes_{{\mathcal O}{}}{\mathscr F},$$
satisfying the Leibniz's rule for the product of {\it holomorphic} functions on $X$ and sections of ${\mathscr F}$.

 All the formulas obtained in Section \ref{S:YMandF} from (\ref{Def:Connection}) to (\ref{FormulaProp})  hold for a holomorphic connection on the coherent sheaf ${\mathscr F}$ over the complex analytic manifold $X$, if we substitute ${\mathscr A}$ by ${\mathscr O}$,   substitute ${\mathscr A}^*$ by the corresponding sheaf of holomorphic forms, and assume that $d$ is the holomorphic exterior differential.

\smallskip

\noindent
{\sc The Atiyah class of  a skyscraper sheaf.} Let $p$ be a point of $X$. On $P:=\{ p\}$ we consider the sheaf ${\mathscr O}_P$ defined by ${\mathscr O}_P(P)={\mathbb C}$. Denoting by $i:P\to X$  the inclusion,  the  direct image sheaf $i_*{\mathscr O}_P$ is an ${\mathscr O}$-module, where the  morphism  ${\mathscr O}\to i_*{\mathscr O}_P$ is given by the evaluation at the point $p$. 
 ${\mathscr S}:=i_*{\mathscr O}_P$ is the skyscraper sheaf on $X$ at $p$ with stalk ${\mathbb C}$. The exact sequence
$$0\to{\mathscr I}_{p|X}\to {\mathscr O}\to {\mathscr S}\to 0,$$
with ${\mathscr I}_{p|X}$ the ideal sheaf defining $p$, shows that    ${\mathscr S}$ is coherent sheaf.

\smallskip

\noindent
{\it Example 1.} Let ${\mathscr S}$ be the skyscraper sheaf on $X={\mathbb C}P^1=\{[x_0:x_1]\,|\, x_i\in{\mathbb C}\}$ at the point  $p=[1:0]$. One has the following exact sequence of ${\mathscr O}$-modules
\begin{equation}\label{exactCP1}
0\to{\mathscr O}(-1)\overset{f}{\rightarrow}{\mathscr O}\overset{q}{\rightarrow}{\mathscr S}\to 0,
\end{equation}
where the morphism $f$ is the multiplication by $x_1$ and $q$ is the evaluation map at the point $p$.

For any   ${\mathscr O}$-module ${\mathscr G}$, ${\rm Ext}^i({\mathscr O},\,{\mathscr G})=H^i(X, {\mathscr G})$ \cite[page 234]{Hartshorne}.
 On the other hand, $\Omega^1({\mathscr S})$ is supported  at $p$, 
since $\Omega^1({\mathscr S})_x\simeq\Omega^1_x\otimes_{{\mathscr O}_x}{\mathscr S}_x$ \cite[page 88]{Kas-Sch}, thus, $H^i(X,\,\Omega^1({\mathscr S}))=0$, for $i>0$. Hence, the ext exact sequence associated to (\ref{exactCP1}) is 
$$\to{\rm Hom} ({\mathscr O},\,\Omega^1({\mathscr S}))\overset{\phi}{\rightarrow} 
{\rm Hom} ({\mathscr O}(-1),\,\Omega^1({\mathscr S}))\to
{\rm Ext}^1 ({\mathscr S},\,\Omega^1({\mathscr S}))\to 0.$$
The morphism $\phi$ is defined as folows: Given $\eta\in{\rm Hom}({\mathscr O},\,\Omega^1({\mathscr S}))$, then $\phi(\eta)(\tau)=\eta(x_1\tau)$, for any local section $\tau$ of
 ${\mathscr O}(-1)$. As $x_1\tau$ vanishes at $p$,  then $\phi=0$ and 
 ${\rm Ext}^1 ({\mathscr S},\,\Omega^1({\mathscr S}))\simeq {\rm Hom} ({\mathscr O}(-1),\,\Omega^1({\mathscr S})).$

On the other hand, the skyscraper sheaf ${\mathscr S}$ does not admit holomorphic connections. 
In fact,     let $h$ and $\tilde h$  be  holomorphic functions defined on an open subset $U\subset{\mathbb C}P^1$,
which contains the point    $p=[1:0]$, such that they coincide at the point $p$, $h(p)=\tilde h(p)$, 
but $\partial h(p)\ne\partial\tilde h(p)$. 
 From the definition of the ${\mathscr O}$-module structure of ${\mathscr S}$, it follows that $h\sigma=\tilde h\sigma$ for any section $\sigma\in {\mathscr S}(U)$. Suppose there was such a connection; 
from the Leibniz's rule applied to $h\sigma$ and $\tilde h\sigma$, we would deduce that $\partial h\otimes \sigma=\partial \tilde h\otimes \sigma$ on $U$, and this equality 
  is in contradiction with the fact that the differentials do not coincide at the point $p$.
 That is, the Atiyah class $a({\mathscr S})\ne 0$. 

\smallskip

The argument above is applicable to skyscraper sheaves on more general manifolds and from it,  one  deduces the following proposition.
 \begin{Prop}\label{P:Skyscraper}
 The skyscraper sheaves over a complex analytic manifold do not admit  holomorphic connections. 
\end{Prop}


\subsection{Holomorphic Yang-Mills fields.} \label{Ss:HoloYM}
Given a coherent sheaf ${\mathscr G}$ on $X$, the fibre of ${\mathscr G}$ at $x\in X$ will be denoted ${\mathscr G}(x);$ that is,
 ${\mathscr G}(x)={\mathscr G}_x/{\mathfrak m}_x{\mathscr G}_x,$ where ${\mathfrak m}_x$ is the maximal ideal of ${\mathscr O}_x$. 
If $\tau$ is a section of ${\mathscr G} $, we set $\tau(x)$ for the image of $\tau_x$ in the fibre ${\mathscr G}(x)$. A Hermitian metric on ${\mathscr G}$ is a family of Hermitian metrics $\langle\,\cdot\;,\,\cdot\,\rangle_x$ on the vector spaces  ${\mathscr G}(x)$, such that for $\tau$ and $\tau'$ sections of ${\mathscr G}$ on an open set $U$, the map $x\mapsto \langle\tau(x),\,\tau'(x)\rangle_x\in{\mathbb C}$ is $C^{\infty} $ on $U\setminus Y$, with $Y$ the singularity set of ${\mathscr G}$.

Let $\nabla$ be a holomorphic connection on  the sheaf ${\mathscr F}$. Any other holomorphic connection is of the form
$\nabla+{\mathcal E}$, with 
\begin{equation}\label{mathcalB}
{\mathcal E}\in{\rm Hom}_{\mathscr O}\big({\mathscr F},\,\Omega^1\otimes_{\mathscr O}{\mathscr F}\big)=
\Gamma\big(X,\,{\mathscr Hom}({\mathscr F},\,\Omega^1\otimes_{\mathscr O}{\mathscr F})\big).
\end{equation}
The singularity set $Z$ of ${\mathscr F}$ is a closed analytic subset of $X$ with ${\rm codim}\,Z\geq 1$ \cite{Scheja}. Moreover, ${\mathscr F}|_{X\setminus Z}$ is a locally free ${\mathscr O}_{X\setminus  Z}$-module.
Thus, on $X\setminus Z$,
${\mathcal E}$ determines a unique  element $E\in \Gamma\big( X\setminus Z,\, \Omega^1\otimes_{\mathscr O}{\mathscr End}(\mathscr F)  \big)$
 (see Lemma  \ref{L:1}).  We will denote by $\psi$ the map 
$$\psi: {\rm Hom}_{\mathscr O}\big({\mathscr F},\,\Omega^1\otimes_{\mathscr O}{\mathscr F}\big)\to \Gamma\big( X\setminus Z,\, \Omega^1\otimes_{\mathscr O}{\mathscr End}(\mathscr F)  \big),\;\;   {\mathcal E}\mapsto E.$$
 
Similarly, the curvature ${\mathcal K}_{\nabla}\in {\rm Hom}_{\mathscr O}\big({\mathscr F},\,\Omega^2\otimes_{\mathscr O}{\mathscr F}\big)$ and its restriction to $X\setminus Z$ defines a vector $K_{\nabla}$ of the vector space $\Gamma\big( X\setminus Z,\, \Omega^2\otimes_{\mathscr O}{\mathscr End}(\mathscr F)  \big)$, (as $X\setminus Z$ is not compact, that vector space may be infinite dimensional).
On $X\setminus Z$ the curvature of $\nabla_{\epsilon}:=\nabla+\epsilon{\mathcal E}$ is 
$K_{\nabla+\epsilon{\mathcal E}}=K_{\nabla}+\epsilon\nabla(E)+O(\epsilon^2).$

 When $X$ is a compact K\"ahler manifold and ${\mathscr F}$ is a holomorphic Hermitian sheaf,  we define
$$||K_{\nabla_{\epsilon}}||^2:=\int_{X\setminus Z}\langle K_{\nabla_{\epsilon}},\, K_{\nabla_{\epsilon}}\rangle\,{\rm d vol}.$$
If $\nabla$ and $\nabla_{\epsilon}$ are compatible with the Hermitian structure, then 
$E\in \Gamma\big( X\setminus Z,\, \Omega^1\otimes_{\mathscr O}\mathfrak{hu}(\mathscr F)  \big).$ We say that $\nabla$ is a {\em holomorphic Yang-Mills field} on ${\mathscr F}$ if 
$$\frac{d}{d\epsilon}\Big|_{\epsilon=0}||K_{\nabla_{\epsilon}}||^2=0,\;\; \text{for all}\;\; E\in {\rm im}(\psi)\cap \Gamma\big( X\setminus Z,\, \Omega^1\otimes_{\mathscr O}\mathfrak{hu}(\mathscr F)  \big).$$
That is, $\nabla$ is a Yang-Mills connection if
\begin{equation}\label{YanMillsEq}
({ K}_{\nabla},\, \nabla{ E}):=\int_{X\setminus Z}\langle{ K}_{\nabla},\, \nabla { E}\rangle\, {\rm d vol}=0,
\end{equation}
where
$\nabla { E}$ is the covariant derivative of ${ E}$.




\smallskip
\noindent
 {\sc Yang-Mills fields on reflexive sheaves.} Let ${\mathscr H}$ be  a coherent reflexive sheaf 
on the complex manifold $X$.  Some   properties satisfied by ${\mathscr H}$ that we will use  are the following. 
 \begin{enumerate}
\item 
 The singularity set $Z$ of ${\mathscr H}$ has codimension grater than $2$ \cite[Cor. 1.4]{Hartshorne1}, and the restriction 
${\mathscr H}|_{X\setminus Z}$ is a locally free sheaf.
\item If $Y\subset X$ is a closed subset of $X$ of codimension $\geq 2$, then the restriction map 
$\Gamma(X,\,{\mathscr H})\to \Gamma(X\setminus Y,\,{\mathscr H})$
 is an isomorphism \cite[Prop. 1.11]{Hartshorne2}. 
 \item 
 If ${\mathscr G}$ is a coherent sheaf, then ${\mathscr Hom}({\mathscr G},\,{\mathscr H})$ is also reflexive \cite[Chap V, Prop. 4.15]{Koba}.
\end{enumerate}

On the other hand, a coherent sheaf ${\mathscr H}$ is reflexive iff there exists an exact sequence
$$0\to{\mathscr H}\to {\mathscr L}\to {\mathscr T},$$ 
with ${\mathscr L}$ locally free and ${\mathscr T}$ torsion-free \cite[Prop. 1.1]{Hartshorne1}.
The functor $\Omega^k\otimes_{\mathscr O}-$ is exact, since $\Omega^k$ is locally free. Thus, if ${\mathscr H}$ is reflexive, we have the exact sequence
$$0\to\Omega^k\otimes_{\mathscr O}{\mathscr H}\to \Omega^k\otimes_{\mathscr O}{\mathscr L}\to\Omega^k\otimes_{\mathscr O}  {\mathscr T}.$$
As $\Omega^k\otimes_{\mathscr O}{\mathscr L}$ is locally free and $\Omega^k\otimes_{\mathscr O}  {\mathscr T}$  torsion-free, it follows that 
$\Omega^k\otimes_{\mathscr O}{\mathscr H}$ is also reflexive.

Let ${\mathscr F}$ be a reflexive sheaf on $X$, whose singularity set is denoted by $Z$, and ${\mathcal C}\in{\rm  Hom}({\mathscr F},\,\Omega^k\otimes_{\mathscr O}{\mathscr F})\big)$. 
From   (2) together with  (3), it follows that ${\mathcal C}$ is determined by its restriction in 
$\Gamma\big(X \setminus Z,\,{\mathscr Hom}({\mathscr F},\,\Omega^k\otimes_{\mathscr O}{\mathscr F})\big)$.
 According to (1)
 $$\Gamma\big(X \setminus Z,\,{\mathscr Hom}({\mathscr F},\,\Omega^k\otimes_{\mathscr O}{\mathscr F})\big)\simeq\Gamma\big( X\setminus Z,\, \Omega^k\otimes_{\mathscr O}{\mathscr End}(\mathscr F)  \big).$$
 Again (2) and (3) give rise to the isomorphism
$$ \Gamma\big( X\setminus Z,\, \Omega^k\otimes_{\mathscr O}{\mathscr End}(\mathscr F)  \big)\simeq  \Gamma\big( X,\, \Omega^k\otimes_{\mathscr O}{\mathscr End}(\mathscr F)  \big).$$
Hence, ${\mathcal C}$ determines and element $C$ of the vector space $\Gamma\big( X,\, \Omega^1\otimes_{\mathscr O}{\mathscr End}(\mathscr F)  \big).$ That is,
\begin{Lem}\label{Lem:finite}
If ${\mathscr F}$ is a reflexive sheaf, then  the map 
$${\mathcal C}\in {\rm Hom}_{\mathscr O}\big({\mathscr F},\,\Omega^k\otimes_{\mathscr O}{\mathscr F}\big)\mapsto
C\in \Gamma\big( X,\, \Omega^k\otimes_{\mathscr O}{\mathscr End}(\mathscr F)  \big)$$
is a bijective correspondence.
\end{Lem}

 Let $\nabla$ be a holomorphic connection on the reflexive sheaf ${\mathscr F}$, whose singularity set is denoted by $Z$.
As we said, any other holomorphic connection on ${\mathscr F} $ is of the form 
$\nabla+{\mathcal E}$, where ${\mathcal E}$ satisfies (\ref{mathcalB}). By Lemma \ref{Lem:finite}, ${\mathcal E}$ determines
 $E\in \Gamma\big( X,\, \Omega^1\otimes_{\mathscr O}{\mathscr End}(\mathscr F)  \big).$
Similarly, the curvature ${\mathcal K}_{\nabla}$ is determined by the vector $K_{\nabla}$ of the  vector space 
$\Gamma\big( X,\, \Omega^2\otimes_{\mathscr O}{\mathscr End}(\mathscr F)  \big)$.



 Let us assume that ${\mathscr F}$ is reflexive sheaf equipped with a Hermitian metric,  and $\nabla$ a Hermitian connection. 
 By the   properties above mentioned, 
  $$\Gamma(X\setminus{Z},\,\Omega^k\otimes_{\mathscr O}\mathfrak{hu}({ F}))\subset  \Gamma\big( X,\setminus Z, \Omega^k\otimes_{\mathscr O}{\mathscr End}(\mathscr F)  \big)=\Gamma\big( X,\, \Omega^k\otimes_{\mathscr O}{\mathscr End}(\mathscr F)  \big),$$ and these spaces are
	 {\it finite dimensional}, since $\Omega^k\otimes_{\mathscr O}{\mathscr End}(\mathscr F)$ is coherent. We will denote  $r_k={\rm dim}\,H^0(X\setminus {Z},\; \Omega^k\otimes_{\mathscr O}\mathfrak{hu}({ F}))$, thus
 $ { K}_{\nabla}$ is a vector of a vector space with dimension $r_2$. 

If  $\nabla$ is a Yang-Mills field on the reflexive sheaf ${\mathscr F} $, then $K_{\nabla}$
		is a vector of an $r_2$-dimensional Hermitian vector space which satisfies the   orthogonality conditions (\ref{YanMillsEq}). When $r_1\geq r_2$, we say that $\nabla$ is a {\it generic} gauge field if the linear map 
	\begin{equation}\label{nablaextendida}
	\nabla:\Gamma(X\setminus{Z},\,\Omega^1\otimes_{\mathscr O}\mathfrak{hu}({ F}))\to \Gamma(X\setminus{Z},\,\Omega^2\otimes_{\mathscr O}\mathfrak{hu}({ F}))
	\end{equation}
	is surjective.
	Thus, we have the following proposition that generalizes Proposition \ref{P:generic}.
\begin{Prop}\label{reflexive}
If ${\mathscr F}$ is a reflexive sheaf endowed with a Hermitian metric and $\nabla$ is a generic Yang-Mills connection, then $\nabla$ is flat.
\end{Prop}

If a coherent sheaf is endowed with a holomorphic flat connection, then it is really a locally free sheaf \cite[Prop. 2.3, page 184]{Gelfand-Manin2}
 \cite[Prop. 2.2.5]{Hotta}. Hence, if $r_1\geq r_2$ and ${\mathscr F}$ is reflexive with singularities, then it does not admit a generic Yang-Mills gauge field. In other words,  if $r_1\geq r_2$ and ${\mathscr F}$ is reflexive with singularities, then the kernel of the morphism (\ref{nablaextendida}) has dimension $>r_1-r_2$, for any Yang-Mills connection $\nabla$. 


\subsection{Meromorphic gauge fields and geometric phases}\label{meromorphicgauge}
Let ${\mathscr F}$ be a coherent ${\mathscr O}{}$-module on the complex analytic manifold $X$.
 We    denote by $Z$ the singularity space of ${\mathscr F}$ and we set  ${\mathscr O}{}[*Z]$ for the sheaf of meromorphic functions on $X$, that are holomorphic on $X\setminus Z$ and have poles 
on $Z$. We set   $\Omega^1{}[*Z]$ for 
$\Omega{}^1\otimes_{{\mathscr O}{}}{\mathscr O}{}[*Z]$ and ${\mathscr F}^{Z}:={\mathscr O}{}[*Z]\otimes_{{\mathscr O}{}}{\mathscr F}$.
A connection  on ${\mathscr F}$ {\it meromorphic} along $Z$ is a morphism of ${\mathbb C}_{X}$-modules
\begin{equation}\label{MeromorphicConn}
\nabla:{\mathscr F}^Z\to\Omega^1{}\otimes _{{\mathscr O}{}}{\mathscr F}^Z,
\end{equation}
which satisfies the Leibniz's rule for the product of functions of ${\mathscr O}{}[*Z]$  and sections of ${\mathscr F}^Z$.

Since the restriction  of the meromorphic connection (\ref{MeromorphicConn}) to the locally free sheaf ${\mathscr F}|_{X\setminus Z}$ is a holomorphic one, the vanishing of the Atiyah class $a({\mathscr F}|_{X\setminus Z})$ is a necessary condition for the existence of  such a meromorphic connection. 

A meromorphic connection whose ringularities are not ``wild'' is called regular; that is, when the equation for the parallel transport is regular in the Fuchs' sense
is regular. In this case the solutions to this equation have a ``moderate'' growth. The  regular flat connections  can be characterised from a topological point of view, as explained in the following paragraph.

\smallskip

\noindent
{\it Flat connections.} The meromorphic connection $\nabla$ is called flat if for any holomorphic vector fields $v,v'$ on $X$ 
$$[\nabla_ v,\,\nabla_{v'}]=\nabla_{[v,\,v']}.$$

 When $Z$ is a divisor of $X$, the flat connection defines a $D_X$-module structure on
${\mathscr F}^Z$ \cite[page 140]{Hotta}.  A regular flat connection   determines a representation of the fundamental group $\pi_1(X\setminus Z)$, according to the Deligne's version of the Riemann-Hilbert correspondence \cite[page 149]{Hotta}. That representation is defined, by means of the monodromy of parallel transport. We will consider that property of the regular flat connections in a brief analysis of the Aharonov-Bohm effect.

\subsubsection{Aharonov-Bohm effect}
The Aharonov-Bohm effect is the change in the phase of the wave function of an electron,  
when it travels in a closed curve around a solenoid ${\it s}$.  The magnetic field   created by ${\it s}$  is confined within it, thus the  potential vector is a {\it closed} $1$-form outside ${\it s}$. 

 Let us consider an ideal solenoid 
indefinitely long with infinitesimal radius, and an electron ${\it e}$, moving in a plane $\pi$ orthogonal to the solenoid.   Even though the confinement of the magnetic field, 
after completing a closed curve around the solenoid,   appears a shift in the phase of the wave function of ${\it e}$. The corresponding factor of phase is ${\rm e}^{-i\Phi}$, where $\Phi$ is the magnetic flux through 
the solenoid \cite[Sec. 4]{Vina20}. 

 A way to deduce this result, in a geometric setting, consists of endowing the plane $\pi$ with a complex structure, and then to express the potential vector for the magnetic field as    a {\it meromorphic} $1$-form on $\pi$. The shift of phase is the integral of this form along the corresponding closed curve.

This effect can also be analyzed in the frame of the meromorphic gauge fields on coherent sheaf as follow.
On $X={\mathbb C}$ we consider the sheaf ${\mathscr G}={\mathscr O}[*Z]$, with $Z=\{0\}$.
 By means of  the meromorphic closed $1$-form $\beta=\frac{-k}{z}dz$, with $k$ a constant in ${\mathbb C}\setminus {\mathbb Z}$, one may define the meromorphic integrable connection on ${\mathscr G}$
$$\nabla(h)=dh+h\otimes\beta,$$
 for $h\in {\mathscr O}{}[*Z]$. 

In this case, the equation for the parallel transport is $dh+h\beta=0$, and it has not a holomorphic solution.  Its solution,  up to multiplicative constant, is
$\tilde h(z)={\rm exp}(k\log z)$, where $\log z$ is any branch of the logarithm function. 
 Let us consider the complex of ${\mathbb C}{}$-modules
$${\mathcal D}^{\bullet}: 0\to{\mathcal D}^{-1}:={\mathscr G}\overset{\nabla}{\longrightarrow}{\mathcal D}^{0}:=
\Omega^1\otimes_{\mathscr O}{\mathscr G}\to 0.$$
This is in fact, the de Rham complex defined by the flat connection $\nabla$ \cite[page 103]{Hotta}.

The cohomology sheaf ${\mathscr H}^{0}({\mathcal D}^{\bullet})$ is trivial, and ${\mathscr H}^{-1}({\mathcal D}^{\bullet})$  is the extension by zero to $\{0\}$ 
of the local system on ${\mathbb C}\setminus\{0\}$ determined by the multivalued function $\tilde h$. 

The analytic continuation of $\tilde h(z)$ along the curve $C=\{{\rm e}^{it}\}_{t\in[0,\,2\pi]}$ changes $\tilde h(z)$ into 
$\tilde h(z)\,{\rm e}^{2\pi k i}$. That is, the above mentioned  local system is the one defined by the representation of $\pi_1(X\setminus \{0\})$ that assigns to the homotopy class $[C]$  the complex number ${\rm e} ^{2\pi k i}$. Which is the factor of phase that appears in the Aharonov-Bohm effect, when $k=\frac{-\Phi}{2\pi}$.

The crucial points are that $\pi_1({\mathbb C}\setminus\{0\})\ne 1$ and that the holomorphic form $\beta|_{{\mathbb C}\setminus\{0\}}$ can not be extended to a holomorphic form on ${\mathbb C}$.
 Hence, the above result does not always  admit of a direct generalization to 
dimensions greater than $1$. This question is analyzed in the following example.

\smallskip
\noindent
{\it Example 2.} Let $Z$ be a submanifold  of the $n$-dimensional complex simply connected manifold $X$, such that  ${\rm codim }\,Z\geq 2$. Thus, $\pi_1(X\setminus Z)=1$.
 We assume also that 
 $X $ is $(n-1)$-complete 
\cite[page 235]{A-G}. In particular, $X$ is cohomologically $(n-1)$-complete; i.e. $H^j(X,\,{\mathscr G} )=0$ for $j\geq n-1$ and for any coherent sheaf ${\mathscr G}$ \cite[page 250]{A-G}. Under these hypotheses, the Hatogs' extension theorem holds on $X$ \cite{C-R, M-P}.

 We set ${\mathscr B}:={\mathscr O}[*Z]$.
As above, let us consider a meromorphic flat connection on  ${\mathscr B}$   defined by a closed $1$-form $\alpha=f\gamma$, with $f\in{\mathscr O}[*Z]$ and $\gamma$ a holomorphic $1$-form. As $f$ is holomorphic on $X\setminus Z$, by the Hatogs' extension theorem, $f$ is   holomorphic on $X$ and the monodromy of the parallel transport equation $dh+h\alpha=0$ is trivial.
Therefore, an electron moving  moving on $X$ in presence of a  magnetic field confined in the submanifold $Z$ will not undergo any change in its wave function.
That is, in this case, there is no ``Aharonov-Bohm'' effect.


\subsubsection{Aharonov-Casher effect}
The Aharonov-Casher effect admits a similar interpretation. This effect appears when a nonrelativistic neutron ${\it n}$ completes a closed curve around an 
electrically charged conductor.  The conductor is supposed to be ideal; that is, unlimited long, infinitesimally thin and with a uniform and static distribution of electric charge. Thus, we can suppose that the conductor is on the $z$-axis of a coordinate system.   When the neutron moves on  a plane orthogonal to the conductor,
describing a closed curve around the wire, the corresponding factor of phase in the spin of ${\it n} $ is ${\rm exp}(i\pi\Lambda \sigma_3)$, where $\Lambda$ is a constant related with the density of charge in the conductor and $\sigma_3$ is the corresponding Pauli 
matrix \cite{Vina20, C-J}. 

In this case, as potential can be taken   a closed {\it meromorphic} $\mathfrak{sl}(3)$-valued $1$-form on the mentioned plane. 
The phase shift can be obtained by integration this form along the curve.


\subsubsection{Wong equation.} 
 The Wong equations describe the variation of the spin of a particle moving in a gauge field \cite{Wong}. Let $p$ be a particle carrying a   like spin variable $I$ with values in the semisimple Lie  algebra ${\mathfrak g}$. Let $\{\tau_a\}_{a=1,\dots, r}$ be a basis of ${\mathfrak g}$ and
$I=\sum_a I^a\tau_a$, then the variation of $I$ when $p$ is moving in the ${\mathfrak g} $-valued gauge potential $A=\sum A^a\tau_a$ satisfies
$$\Dot{ I}^a+\sum C^a_{bc}A^b(\Dot x)I^c=0,$$
where $x(t)$ is the curve described by $p$, and the $C^a_{bc}$ are the structure constants of ${\mathfrak g} $ relative to the basis 
$\{\tau_a\}$ \cite[p. 53]{Bala}, \cite{Magazev}.

The above equations can be written as
\begin{equation}\label{Wong0}
\frac{dI}{dt}+ [A(\Dot x(t)),\, I(t)]=0,
\end{equation}
where $[\,\cdot\, ,\,\cdot\,]$ is the Lie bracket in ${\mathfrak g}$.

Since ${\mathfrak g}$ is semisimple, the adjoint representation ${\rm ad}:{\mathfrak g}\to{\rm End}({\mathfrak g})$ is a monomorphism of algebras. Thus, we   will {\it  identify} a vector $u\in {\mathfrak g}$ with the linear map ${\rm ad}(u)$. Applying ${\rm ad}$ to (\ref{Wong0}), one obtains the following equation 
$$\frac{dI}{dt}+ A(\Dot x(t))\circ I(t)-I(t)\circ A(\Dot x(t))=0,$$
which can be regarded as the parallel transport of the endomorphism $I$ along the curve $x(t)$ (see (\ref{nablaEnd}) for the case $k=0$). If the gauge field $A$ is {\it singular}, then it is to be expected phase shifts in the variable $I$, when $p$ completes a closed curve. 

A geometric setting for analyzing Wong equation is the following.
 We denote by ${V}$ the trivial bundle
 $X\times {\mathfrak g}_{\mathbb C}$ over the   complex  manifold $X$, and  will regard ${\mathfrak g}$ as a subalgebra of ${\rm End}({\mathfrak g})$. Then, as
 the vectors of ${\mathfrak g}$ define endomorphisms of ${V}$,  
one can consider   connections 
$$\nabla:{V}\to \Omega^1\otimes_{\mathscr O}{\mathfrak g}_{\mathbb C}.$$

 Letting $s_a(x):=(x,\,\tau_a)$, 
 the set
 $\{s_a\}$ is a global frame for ${ V}$. Hence,
$$\nabla(s_i)=\sum_{j=1}^r\alpha_{ji}\otimes s_j,\;\; i=1,\dots, r,$$
 with $\alpha_{ij}$ $1$-forms on $X$. 

 Let us assume that the particle $p$  is moving along the curve $x(t)$ on the manifold $X$ in presence of the gauge field $\nabla$. As we said, the variation of the variable spin like $I$ is determined by the Wong equation
$$\nabla_{\Dot x(t)}I=0,$$
which in turn defines the corresponding holonomy when $x(t)$ is a closed curve. We will consider some particular cases.  
\begin{enumerate}
 
\item
The shift factor in the variable $I$,  after the particle has completely traveled   the curve, is given by the respective holonomy. If the connection $\nabla$ is  holomorphic and flat,  and the loop is null-homotopic, then the shift factor is trivial.

\item
Let us assume that $X$ and $Z$ satisfy the properties stated in Example 2. Moreover, we suppose that
 $\alpha_{ji}=f_{ji}\gamma_{ji}$, with $f_{ji}\in {\mathscr O}[*Z]$ and $\gamma_{ji}$ holomorphic form on $X$. By the Hartogs' extension theorem, $\alpha=(\alpha_{ji})$ is a matrix of holomorphic $1$-forms. If $d\alpha+\alpha\wedge\alpha=0$, then the curvature of the connection vanishes and the    shift factor is the identity.

\item
If $\nabla$ is a meromorphic regular flat connection, in the sense that the $\alpha_{ij}\in\Omega^1[*Z]$, with $Z$ a divisor of $X$, and  the parallel differential equation is regular in Fuchs' sense, 
 then the shift factors are given by a representation of $\pi_1(X\setminus Z)$. This is a consequence of the Deligne's version of the Riemann-Hilbert correspondence \cite[Sect. 5.2.3]{Hotta} \cite{Malgrange}.

\end{enumerate}

\end{document}